\documentclass[12pt,twoside]{article}

\usepackage{amssymb,amsmath,latexsym}
\usepackage{epsfig}
\usepackage{graphicx}

\def\Re{{\rm Re \, }}
\def\Im{{\rm Im \, }}

\DeclareMathOperator*{\Res}{Res}
\DeclareMathOperator{\Ai}{Ai}

\newtheorem{theorem}{Theorem}[section]
\newtheorem{lemma}[theorem]{Lemma}
\newtheorem{corollary}[theorem]{Corollary}
\newtheorem{proposition}[theorem]{Proposition}

\newtheorem{Definition}[theorem]{Definition}
\newenvironment{definition}{\begin{Definition}\rm}{\end{Definition}}
\newtheorem{Remark}[theorem]{Remark}
\newenvironment{remark}{\begin{Remark}\rm}{\end{Remark}}
\newtheorem{Example}[theorem]{Example}

\newenvironment{proof}%
{\rm \trivlist \item[\hskip \labelsep{\bf Proof. }]}%
{\hspace*{\fill}$\Box$\endtrivlist}
{\rm \trivlist \item[\hskip \labelsep{\bf Proof}]}%
{\hspace*{\fill}$\Box$\endtrivlist}

\begin{document}

\begin{center}
 \textbf{\Large
  Riemann--Hilbert analysis for
  Laguerre polynomials with large negative parameter
 }\\[5mm]
 \textit{\large
  A.B.J. Kuijlaars\footnote{Supported in part by FWO research project G.0176.02,
  by INTAS project 00-272, and by a research grant of the
  Fund for Scientific Research--Flanders.}
  and
  K.T-R McLaughlin\footnote{Supported in part by
  NSF grant \# DMS-9970328.}
 }\\[5mm]
\end{center}

\small
\noindent\textbf{Abstract.}
We study the asymptotic behavior of Laguerre polynomials $L_n^{(\alpha_n)}(nz)$
as $n \to \infty$, where $\alpha_n$ is a sequence of negative parameters such
that $-\alpha_n/n$ tends to a limit $A > 1$ as $n \to \infty$.
These polynomials satisfy a non-hermitian orthogonality on certain
contours in the complex plane. This fact allows the formulation of a Riemann--Hilbert
problem whose solution is given in terms of these Laguerre polynomials.
The asymptotic analysis of the Riemann--Hilbert problem
is carried out by the steepest descent method of Deift and Zhou, in the
same spirit as done by Deift et al.\ for the case of orthogonal polynomials
on the real line. A main feature of the present paper is the choice of the
correct contour.

\vspace{2mm}

\vspace{2mm}

\normalsize
\noindent

\section{Introduction}
\subsection{Generalized Laguerre polynomials}

The classical Laguerre polynomials $L_n^{(\alpha)}$ are orthogonal
on the interval $[0,\infty)$ with respect to the weight $x^{\alpha} e^{-x}$,
that is
\begin{equation} \label{eq11}
    \int_0^{\infty} L_n^{(\alpha)}(x) L_m^{(\alpha)}(x)
    x^{\alpha} e^{-x}\, dt = 0, \qquad
    \mbox{ if } n \neq m.
\end{equation}
The integral in (\ref{eq11}) converges only if $\alpha > -1$.
The Laguerre polynomials are given by the explicit formula
\begin{equation} \label{eq12}
    L_n^{(\alpha)}(x)  = \sum_{k=0}^n \binom{n+ \alpha}{n-k}
    \frac{(-x)^k}{k!}
\end{equation}
and by the Rodrigues formula
\begin{equation} \label{eq13}
    L_n^{(\alpha)}(x) = \frac{1}{n!} x^{-\alpha} e^x
        \left(\frac{d}{dx} \right)^n \left[ x^{\alpha + n} e^{-x} \right],
\end{equation}
see e.g.\ \cite{Szego}.
Both (\ref{eq12}) and (\ref{eq13}) make sense for arbitrary $\alpha$ (even complex)
and they define the generalized non-classical Laguerre polynomials.
The recurrence relation
\begin{equation} \label{eq14}
    -x L_n^{(\alpha)}(x) = (n+1) L_{n+1}^{(\alpha)}(x) - (2n+\alpha+1) L_n^{(\alpha)}(x)
       + (n+\alpha) L_{n-1}^{(\alpha)}(x),
\end{equation}
with $L_0^{(\alpha)} \equiv 1$, $L_{-1}^{(\alpha)} \equiv 0$,
and the second order differential equation
\begin{equation} \label{eq15}
    xy''(x) + (\alpha + 1 - x) y'(x) + ny(x) = 0,
    \qquad y(x) = L_n^{(\alpha)}(x),
\end{equation}
continue to hold for arbitrary $\alpha$.  We consider in this paper only real
and negative $\alpha$, although extensions to complex $\alpha$ are possible.

\subsection{Relations to other polynomials}
Laguerre polynomials with negative parameters appear in the literature in a number
of forms. First we note the special cases $\alpha = -n$ and $\alpha = -n -1$,
where we have
\begin{equation} \label{eq16}
    L_n^{(-n)}(z) = \frac{(-1)^n}{n!} z^n,
\end{equation}
and
\begin{equation} \label{eq17}
    L_n^{(-n-1)}(z) = (-1)^n \sum_{k=0}^{n} \frac{z^k}{k!},
\end{equation}
respectively. Thus for $\alpha = -n-1$, the generalized Laguerre polynomial agrees
(up to a sign if $n$ is odd) with the partial sum of the exponential series.
More generally, we have that for $\alpha = - n -m$ with $n,m \in \mathbb N$,
the generalized Laguerre polynomials appear as the numerator and denominator
polynomials in the rational Pad\'e approximant for the exponential function.
To be precise, if
\begin{equation} \label{eq18}
    p_{n,m}(x) = (-1)^n \binom{n+ m}{n}^{-1} L_n^{(-n-m-1)}(x)
\end{equation}
and
\begin{equation} \label{eq19}
    q_{n,m}(x) = (-1)^m \binom{n+m}{m}^{-1} L_m^{(-n-m-1)}(x)
\end{equation}
then $p_{n,m}(0) = q_{n,m}(0) = 1$, and
\begin{equation} \label{eq110}
    p_{n,m}(x) - q_{n,m}(x) e^x = O\left(x^{n+m+1}\right)
    \qquad \mbox{ as } x \to 0,
\end{equation}
see e.g.\ \cite{Perron}, \cite{SV1}.

Laguerre polynomials with negative parameters are also related to the
so-called generalized Bessel polynomials
\begin{equation} \label{eq111}
    y_n(z;a) = \sum_{k=0}^n \binom{n}{k} (n+a-1)_k \left(\frac{z}{2}\right)^k,
\end{equation}
since
\begin{equation} \label{eq112}
    y_n(z;a) = (-1)^n n! \left(\frac{z}{2}\right)^n L_n^{(-2n-a+1)} \left(\frac{2}{z}\right).
\end{equation}
The usual Bessel polynomials correspond to $a=2$ in (\ref{eq111}).
See Grosswald \cite{Gros} for a comprehensive account of these polynomials.

\subsection{Earlier work on asymptotics}
\begin{figure}[th!]
\centerline{\includegraphics[width=5cm]{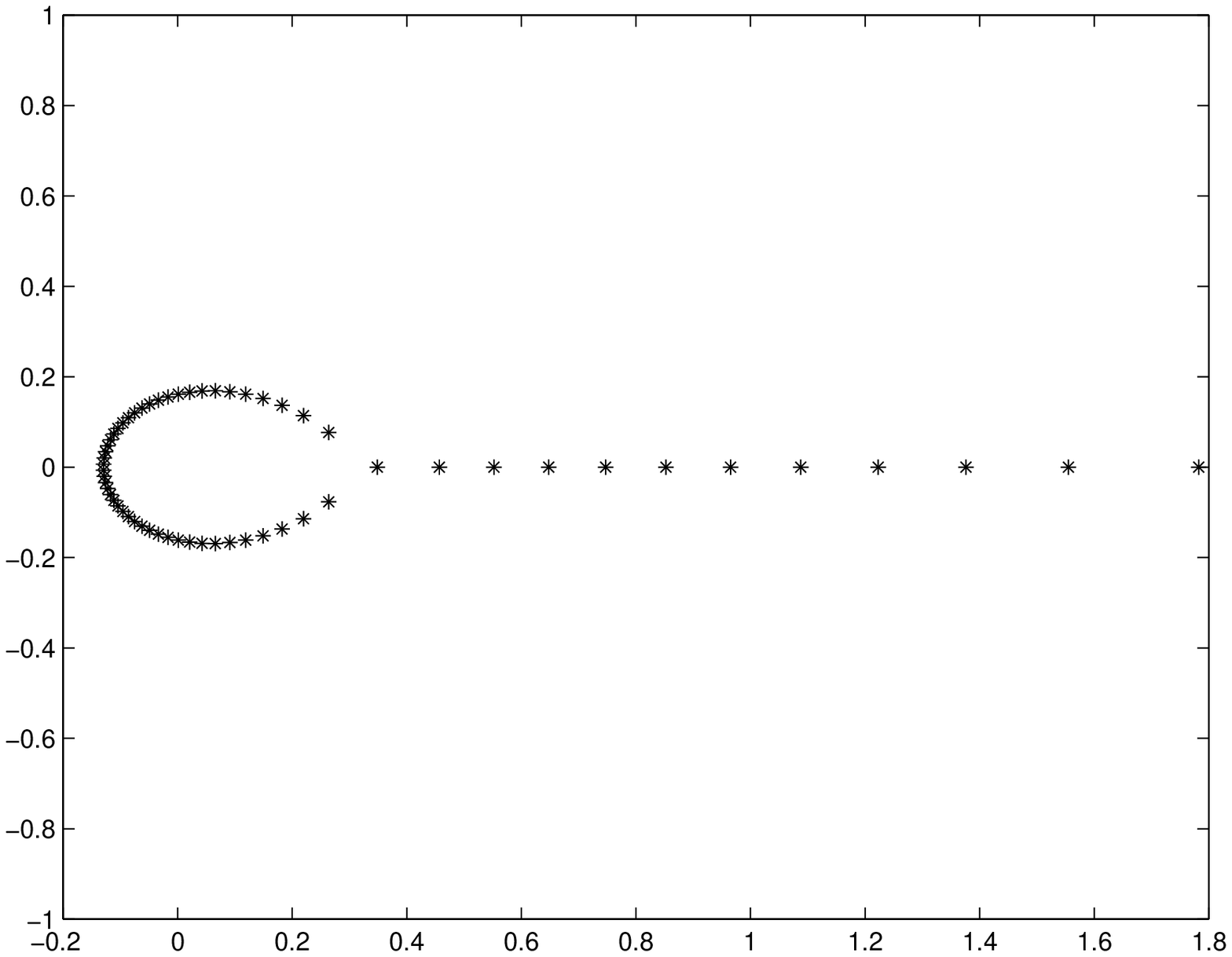}
\includegraphics[width=5cm]{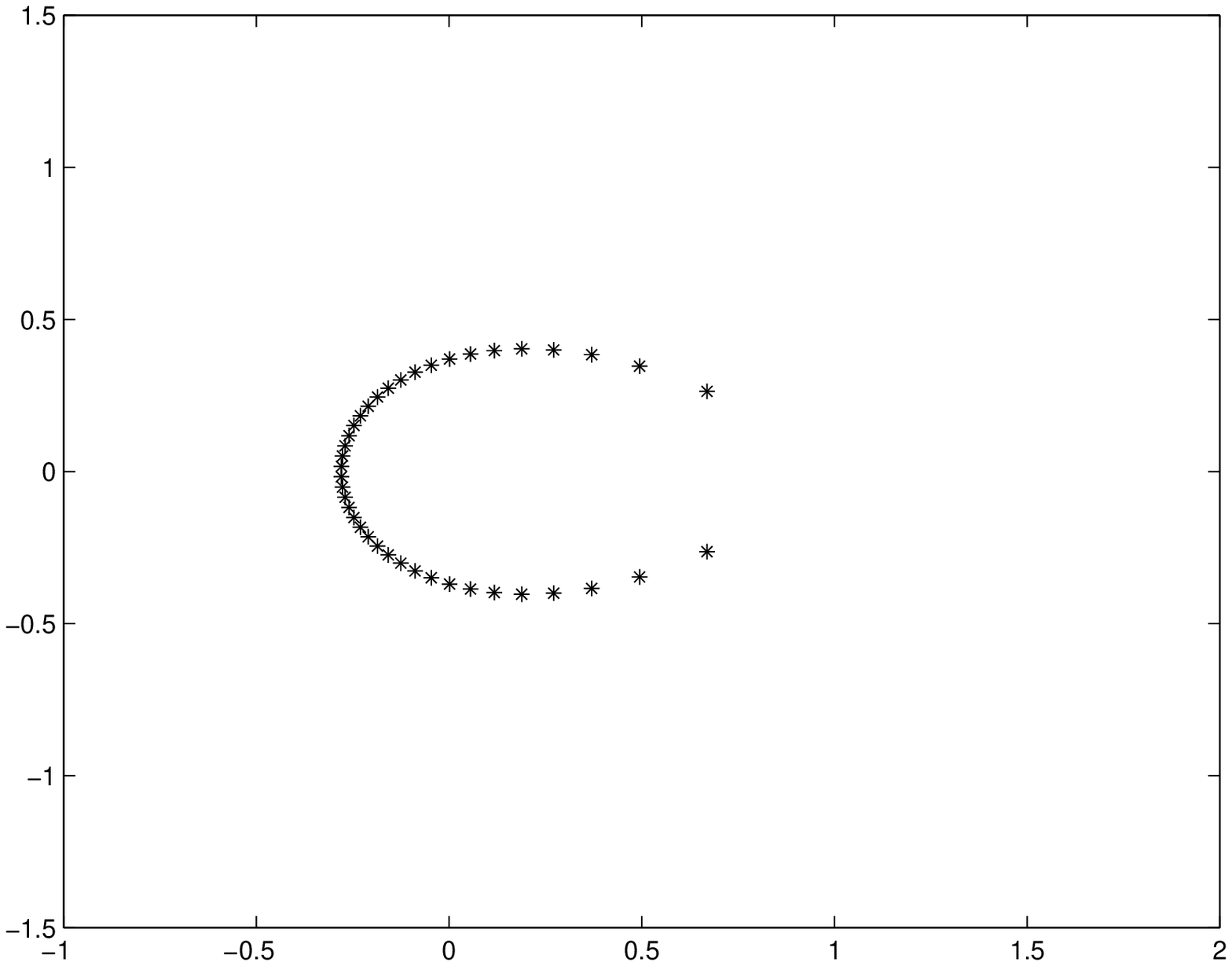}
\includegraphics[width=5cm]{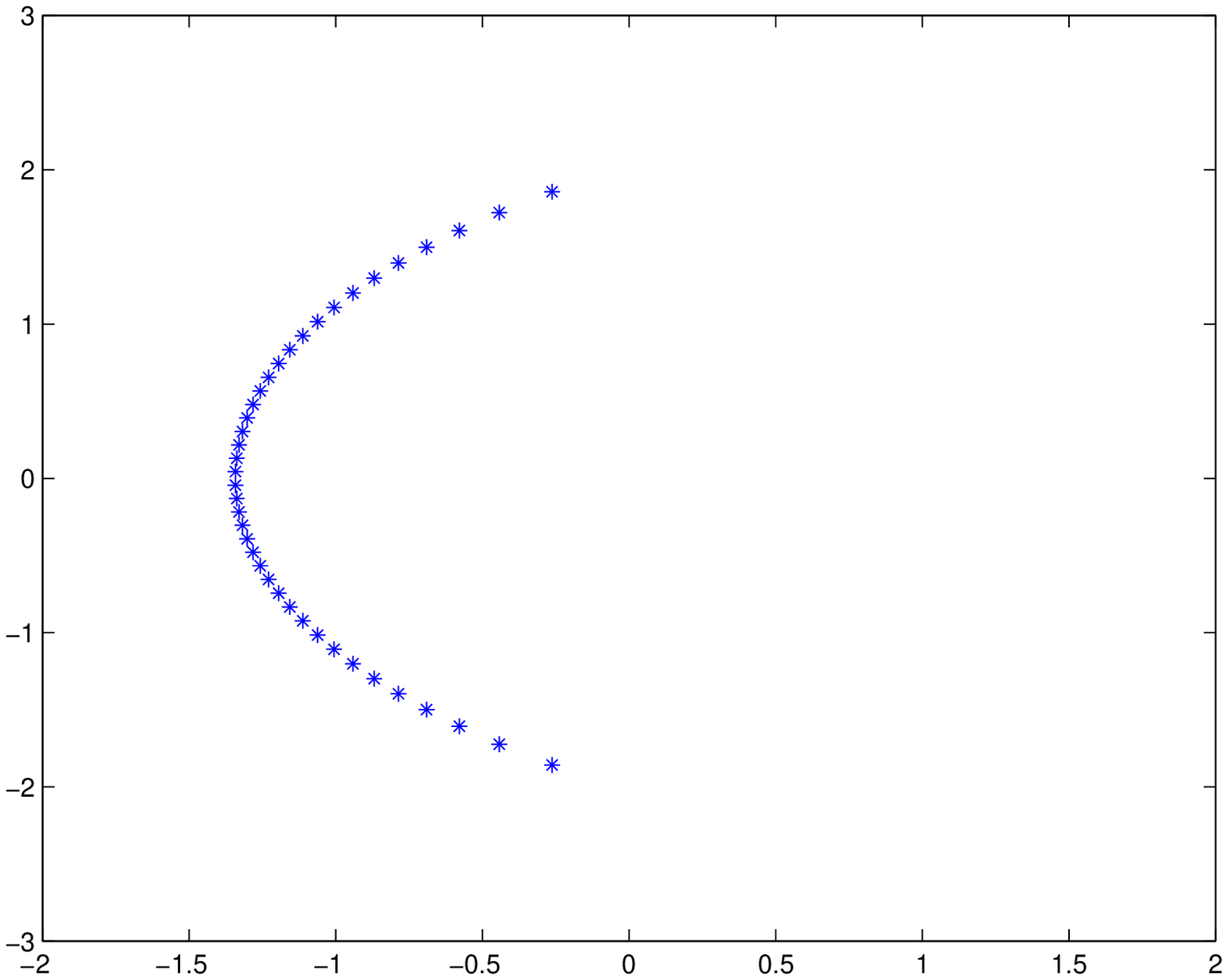}}
\caption{Zeros of generalized Laguerre polynomials $L_{n}^{(-An)}(nx)$ for $n = 40$
and $A = 0.81$ (left), $A = 1.01$ (middle), and $A = 2$ (right).}
\label{fig:plot1}
\end{figure}

For $\alpha > -1$, the Laguerre polynomials satisfy the orthogonality (\ref{eq11})
on the positive real axis, and therefore they have only positive real zeros. This property is
lost for $\alpha < -1$. Indeed, the generalized Laguerre polynomials may
have many non-real zeros. In Figure \ref{fig:plot1}
we have plotted the zeros of $L_{40}^{(-40A)}(40z)$
for a number of values of $A > 0$. Similar plots are shown in the paper
\cite{MMO1}.

In Figure \ref{fig:plot1} we see that the zeros cluster along certain curves
in the complex plane. Mart\'{\i}nez et al.\ \cite{MMO1}
identified these curves as trajectories of a quadratic differential,
depending on a parameter $A$.
For $A > 1$, the curve is a simple arc, which as $A$ decreases to $1$ closes
itself to form for $A = 1$ the well-known Szeg\H{o} curve \cite{Szego1}, \cite{PV}.
For $0 < A < 1$, the curve consists
of a closed loop together with an interval on the positive real axis.

A number of rigorous results on the asymptotic behavior of the zeros
of generalized Laguerre polynomials $L_n^{(\alpha_n)}(nx)$ such that
\begin{equation} \label{eq113}
    \lim_{n \to \infty} \frac{\alpha_n}{n} = -A
\end{equation}
are known from the literature. The first result is due to Szeg\H{o} \cite{Szego1}
who studied the partial sum of the exponential series, that is $\alpha_n = -n-1$
see (\ref{eq17}). Szeg\H{o} showed that the normalized zeros tend to the curve which now
bears his name. Olver \cite{Olver1} considered the zeros of Hankel functions,
which includes the Bessel polynomials as a special case. In terms of the
generalized Laguerre polynomials, this is the case $\alpha_n = -2n-1$.
 Saff and Varga \cite{SV2} studied the zeros and poles of
Pad\'e approximants to the exponential function, see (\ref{eq18})--(\ref{eq19}).
 Their main
result says that for integers $\alpha_n < -n$ such that (\ref{eq113}) holds,
all zeros of $L_n^{(\alpha_n)}(nx)$ tend to a well-defined curve as $n \to \infty$.
The curve depends on $A \geq 1$ only, and coincides with the curve
described in \cite{MMO1}. Saff and Varga also obtained the weak limit of
the zero counting measures. The proofs in \cite{SV2} can be extended without
any difficulty to non-integer $\alpha_n < -n$.

Stated in terms of generalized Bessel polynomials (\ref{eq111})--(\ref{eq112})
asymptotic results on zeros are due to De Bruin et al.\ \cite{dBSV},
Carpenter \cite{Carp}, and Wong and Zhang \cite{WZ}. These papers deal with
the limit (\ref{eq113}) with $A = 2$. The latter paper also presents
uniform asymptotic expansions of the generalized Bessel polynomials.
In very recent work, Dunster \cite{Dunster} establishes uniform
asymptotic expansions in the complex plane for the case of general
$A$, with the exception of $A = 0$ and $A = 1$. The results on zeros
in \cite{Dunster} are restricted to the case $A > 1$.

\subsection{Asymptotics from Riemann--Hilbert problems}
Most papers cited above use some form of the steepest descent technique
for integrals, see especially \cite{SV2} and \cite{WZ}.
Mart\'{\i}nez et al.\ \cite{MMO1} use an orthogonality relation in the
complex plane satisfied by generalized Laguerre polynomials.
The approach of Dunster \cite{Dunster} starts from the
differential equation (\ref{eq15}) and is based on techniques developed
by Olver \cite{Olver2}.

In this paper we derive asymptotics of generalized Laguerre polynomials
using the nonlinear steepest descent / stationary phase method for
Riemann--Hilbert problems introduced by Deift and Zhou in \cite{DZ1},
and further developed in \cite{DZ2} and \cite{DVZ}. In later
developments, the method was applied successfully to problems in random matrix
theory \cite{DIZ}, \cite{DKMVZ1}, and in orthogonal polynomials
\cite{BI}, \cite{DKMVZ1}, \cite{DKMVZ2}, \cite{KM}, and combinatorics
\cite{BDJ}. For review of some of these developments, and a pedagogic
introduction to some of the material of random matrix theory, orthogonal
polynomials, and Riemann--Hilbert problems, see \cite{Deift}.

The Riemann--Hilbert approach to the asymptotics of generalized Laguerre
polynomials starts from the observation
that these polynomials satisfy orthogonality relations in the complex plane.
The orthogonality is on a contour $\Sigma$  going around
the positive real axis, but otherwise being quite arbitrary.
The orthogonality property allows the
formulation of a Riemann--Hilbert problem, due to Fokas, Its,
and Kitaev \cite{FIK}, whose solution is given in terms of $L_n^{(\alpha)}$.
The Riemann--Hilbert problem is analyzed in the large $n$ limit
with the steepest descent method as done in \cite{DKMVZ1} and \cite{DKMVZ2}
for orthogonal polynomials on the real line.

A novel feature for the problem at hand is that the arbitrary
contour $\Sigma$ has to be chosen in a correct way in order to arrive
at a Riemann--Hilbert problem which is amenable to subsequent asymptotic
analysis. The correct contour was described in \cite{MMO1}.
It is a curve with the S-property
of Stahl \cite{Stahl} and Gonchar and Rakhmanov \cite{GR}.
The structure of the curve depends on the value of $A$ as already explained
before.
In this paper we analyze the case of an open contour, that is, the case $A > 1$.
In subsequent work we consider the case of a closed loop plus an interval
($0 < A < 1$) and the case of a single closed contour ($A = 1$).
We note that the steepest descent / stationary phase method for
Riemann--Hilbert problems was augmented to handle cases in which the
contour selection involves determining a set of nontrivial curves in
the plane in \cite{KMM}, by Kamvissis, McLaughlin, and Miller, in the
context of the semi-classical limit of the focusing nonlinear Schr\"odinger
equation.

We emphasize that the main interest in the present paper lies in the method
we use and not in the results obtained for the Laguerre polynomials.
In particular, we do not improve upon the asymptotic expansions of
Dunster \cite{Dunster}. The steepest descent method for Riemann--Hilbert
problems is a very powerful new method, and its use in the study of
classical special functions is new. In future work we consider
generalized Laguerre polynomials for the cases $0 < A < 1$ and $A = 1$,
and the Riemann--Hilbert approach will lead to new results for these
cases.

\section{Complex orthogonality and the formulation of the Riemann--Hilbert problem}
\setcounter{equation}{0}

\subsection{Orthogonality}
For $\alpha < -1$, the generalized Laguerre polynomial $L_n^{(\alpha)}$ is
not orthogonal on the positive real axis, but instead satisfies a
non-hermitian orthogonality in the complex plane.

\begin{figure}[th!]
\centerline{\includegraphics[width=8cm]{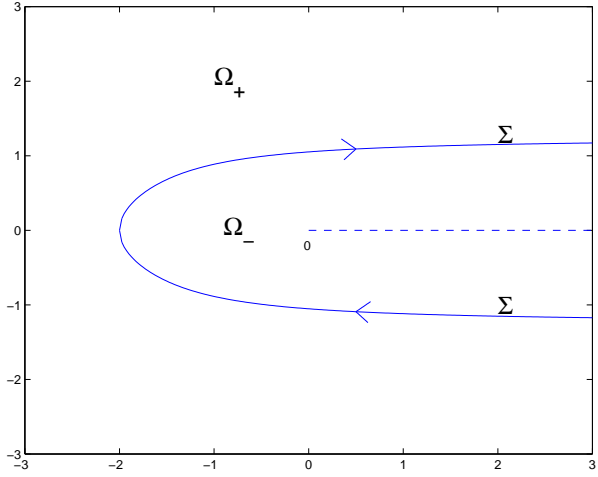}}
\caption{Example of contour $\Sigma$}
\label{fig:sigma}
\end{figure}

Let $\cal F$ be the collection of all simple Jordan curves $\Sigma$ in
$\mathbb C \setminus [0, \infty)$ that are symmetric with respect to the real axis,
and such that there is  $M > 0$ such that for all $x \geq M$, there is
$y(x) > 0$, such that the intersection of $\Sigma$ with $\Re z = x$ consists of the
two points $x \pm i y(x)$, and $\lim_{x \to \infty} y(x) = L$
exists and is finite (possibly $0$).
Any curve $\Sigma \in {\cal F}$ divides the complex plane into two domains,
$\Omega_+$ and $\Omega_-$, where  $\Omega_-$ contains the positive real
axis. We choose the orientation of $\Sigma$ such that $\Omega_+$ is on the $+$-side
(i.e., on the left) while traversing $\Sigma$ and $\Omega_-$ is on the $-$-side.
So $\Sigma$ is oriented clockwise as in Figure \ref{fig:sigma}.

In what follows we define $x^{\alpha}$ with a branch cut along the positive
real axis. Thus  $x^{\alpha} = |x|^{\alpha} e^{i \alpha \arg x}$
with $\arg x \in [0, 2 \pi)$.
\begin{lemma}
Let $\Sigma \in {\cal F}$, $n \in \mathbb N$, and $\alpha \in \mathbb R$.
Then
\begin{equation} \label{eq21}
    \int_{\Sigma} L_n^{(\alpha)}(x) x^k x^{\alpha} e^{-x} \, dx = 0,
    \qquad \mbox{for } k = 0, 1, \ldots, n-1.
\end{equation}
If in addition $\alpha + n + 1 \not\in \mathbb N$, then
\begin{equation} \label{eq22}
     \int_{\Sigma} L_n^{(\alpha)}(x) x^k x^{\alpha} e^{-x} \, dx \neq 0,
     \qquad \mbox{for } k = n.
\end{equation}
\end{lemma}
\begin{proof}
The orthogonality (\ref{eq21}) follows from the Rodrigues formula
(\ref{eq13}) by repeated integration by parts. In the same way, we also get
\[ \int_{\Sigma} L_n^{(\alpha)}(x) x^n x^{\alpha} e^{-x} \, dx =
    \frac{1}{n!} \int_{\Sigma} \left(\frac{d}{dx}\right)^n \left(x^{\alpha + n} e^{-x} \right)
    x^n dx
    = (-1)^n \int_{\Sigma} x^{\alpha+n} e^{-x} \, dx
\]
For $\alpha + n > -1$, we deform $\Sigma$ to the positive real axis to obtain
\begin{eqnarray} \nonumber
\int_{\Sigma} L_n^{(\alpha)}(x) x^n x^{\alpha} e^{-x} \, dx
& = & (-1)^n \left(1- e^{2\pi i \alpha}\right) \int_0^{\infty} x^{\alpha + n} e^{-x} \, dx \\
\label{eq23}
&= &  (-1)^{n+1} 2i e^{\pi i \alpha} \sin (\pi \alpha)  \Gamma(\alpha + n + 1),
\end{eqnarray}
where $\Gamma$ denotes the Gamma function.
By analytic continuation the integral in (\ref{eq22}) is equal to (\ref{eq23})
for every $\alpha$, and (\ref{eq22}) follows.
\end{proof}
The formula (\ref{eq21}) expresses orthogonality with
respect to the complex measure $x^{\alpha} e^{-x}dx$ on $\Sigma$.

\subsection{Riemann--Hilbert problem}
Let $\alpha \in \mathbb R$. We consider the monic polynomials
\begin{equation} \label{eq24}
    P_n(z) = \frac{(-1)^n n!}{n^n} L_n^{(\alpha)}(nz), \qquad n = 0, 1, \ldots,
\end{equation}
Introducing a change of variables $x =nz$ in (\ref{eq21}) and (\ref{eq22}),
we see that
\begin{equation} \label{eq25}
    \int_{\Sigma} P_n(z) z^k z^{\alpha} e^{-nz}\, dz
    \left\{ \begin{array}{ll}
    = 0, & \quad  \mbox{ for } k = 0, 1, \ldots, n-1, \\[10pt]
    \neq 0, & \quad \mbox{ for } k = n,
    \end{array} \right.
\end{equation}
for every contour $\Sigma \in {\cal F}$, provided that
$\alpha +n +1 \not\in \mathbb N$.

The polynomial $P_n$ is characterized through a Riemann--Hilbert problem
due to Fokas, Its, and Kitaev \cite{FIK}.

\subsubsection*{Riemann--Hilbert problem for $Y$:}
Let $\Sigma$ be a contour from the class ${\cal F}$, that divides the
complex plane into two parts $\Omega_+$ and $\Omega_-$, as above.
The problem is to determine a $2 \times 2$
matrix valued function $Y : \mathbb C \setminus \Sigma \to \mathbb C^{2 \times 2}$
such that the following hold.
\begin{enumerate}
\item[(a)] $Y(z)$ is analytic for $z \in \mathbb C \setminus \Sigma$,
\item[(b)] $Y(z)$ possesses continuous boundary values for $z \in \Sigma$,
denoted by $Y_+(z)$ and $Y_-(z)$, where $Y_+(z)$ and $Y_-(z)$ denote
the limiting values of $Y(z')$ as $z'$ approaches $z \in \Sigma$ from $\Omega_+$
and $\Omega_-$, respectively, and
\begin{equation} \label{eq26}
    Y_+(z) = Y_-(z)
    \left( \begin{array}{cc} 1 & z^{\alpha} e^{-nz} \\ 0 & 1 \end{array} \right),
    \qquad \mbox{ for } z \in \Sigma,
\end{equation}
\item[(c)] $Y(z)$ has the following behavior as $z \to \infty$:
\begin{equation} \label{eq27}
    Y(z) = \left( I + O\left(\frac{1}{z}\right) \right)
    \left( \begin{array}{cc} z^n & 0 \\ 0 & z^{-n} \end{array} \right)
    \qquad \mbox{ as } z \to \infty, \ z \in \mathbb C \setminus \Sigma.
\end{equation}
\end{enumerate}

\begin{proposition} Let $\alpha \in \mathbb R$ with $\alpha + n  \not\in \mathbb N$.
Then the unique solution of the Riemann--Hilbert problem for $Y$ is given by
\begin{equation} \label{eq28}
    Y(z) = \left( \begin{array}{ccc}
    P_n(z) & \frac{1}{2\pi i}  \int\limits_{\Sigma}
        \frac{P_n (\zeta) \zeta^{\alpha}e^{-n\zeta}}{\zeta - z} \, d\zeta \\[10pt]
    Q_{n-1}(z) & \frac{1}{2\pi i}  \int\limits_{\Sigma}
     \frac{Q_{n-1}(\zeta) \zeta^{\alpha} e^{-n\zeta}}{\zeta -z} \, d\zeta
     \end{array} \right)
\end{equation}
where $P_n(z)$ is the monic generalized Laguerre polynomial {\rm (\ref{eq24})}
and
\begin{equation} \label{eq29}
    Q_{n-1}(z)= \frac{(-1)^{n+1} n^{n+\alpha} \pi e^{-\pi i\alpha}}
    {\sin (\pi \alpha) \Gamma(\alpha + n)}
    L_{n-1}^{(\alpha)}(nz).
\end{equation}
\end{proposition}
\begin{proof}
The proof is as in \cite[Section 3.2]{Deift}.
Here we will only give the proof for the second row, since that is where
the condition on $\alpha$ plays a role.

From (\ref{eq26}) it follows that $Y_{21}$ is an entire function, which
by (\ref{eq27}) satisfies $Y_{21}(z) = O(z^{n-1})$ as $z \to \infty$. Therefore
$Y_{21}(z) = Q_{n-1}(z)$ for some polynomial $Q_{n-1}$ of degree at most $n-1$.
The (2,2) entry of the condition (\ref{eq26}) is
$(Y_{22})_+(z) = (Y_{22})_-(z) + Q_{n-1}(z) z^{\alpha} e^{-nz}$, which
by the Sokhotskii-Plemelj formula yields
\begin{equation} \label{eq210}
    Y_{22}(z) = \frac{1}{2\pi i}  \int_{\Sigma}
     \frac{Q_{n-1}(\zeta) \zeta^{\alpha} e^{-n\zeta}}{\zeta -z} \, d\zeta.
\end{equation}
From the (2,2) entry of (\ref{eq27}) it follows that $Y_{22}(z) = z^{-n} + O(z^{-n-1})$
as $z \to \infty$. Because of (\ref{eq210}) this gives the conditions
\begin{equation} \label{eq211}
    \int_{\Sigma} Q_{n-1}(\zeta) \zeta^{\alpha} e^{-n\zeta} \zeta^k \, d\zeta = 0
    \qquad \mbox{for } k = 0, \ldots, n-2,
\end{equation}
and
\begin{equation} \label{eq212}
      \int_{\Sigma} Q_{n-1}(\zeta) \zeta^{\alpha} e^{-n\zeta} \zeta^{n-1} \, d\zeta =
    -2\pi i.
\end{equation}
The orthogonality conditions (\ref{eq211}) are satisfied if
$Q_{n-1}(z) = d_n L_{n-1}^{(\alpha)}(nz)$, and the constant $d_n$ must
be chosen so that (\ref{eq212}) holds as well. Because of Lemma 2.1 this
can be done if $\alpha + n \not\in \mathbb N$, and the
result is given by (\ref{eq29}).
\end{proof}

\section{Selection of the right contour}
\setcounter{equation}{0}

\subsection{The right contour $\Sigma$}
In Section 3--6, we assume $n \in \mathbb N$ and $\alpha < -n$ are fixed.
We write $A = -\alpha/n$, so that $A > 1$. [Later, when we let $n \to \infty$,
$\alpha$ and $A$, as well as other notions introduced in these sections,
will depend on $n$.]

\begin{figure}[th!]
\centerline{\includegraphics[width=5cm]{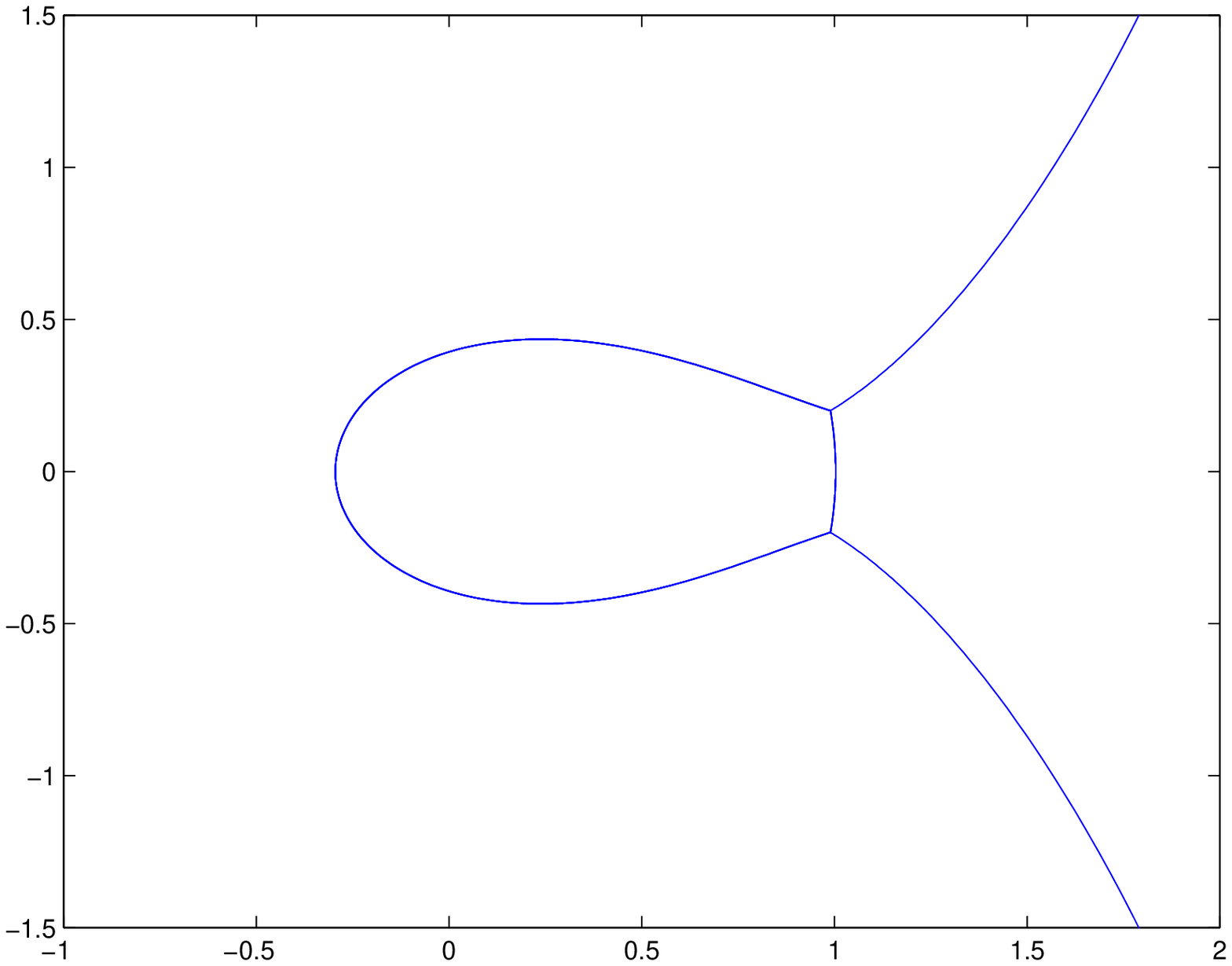} \qquad \qquad
    \includegraphics[width=5cm]{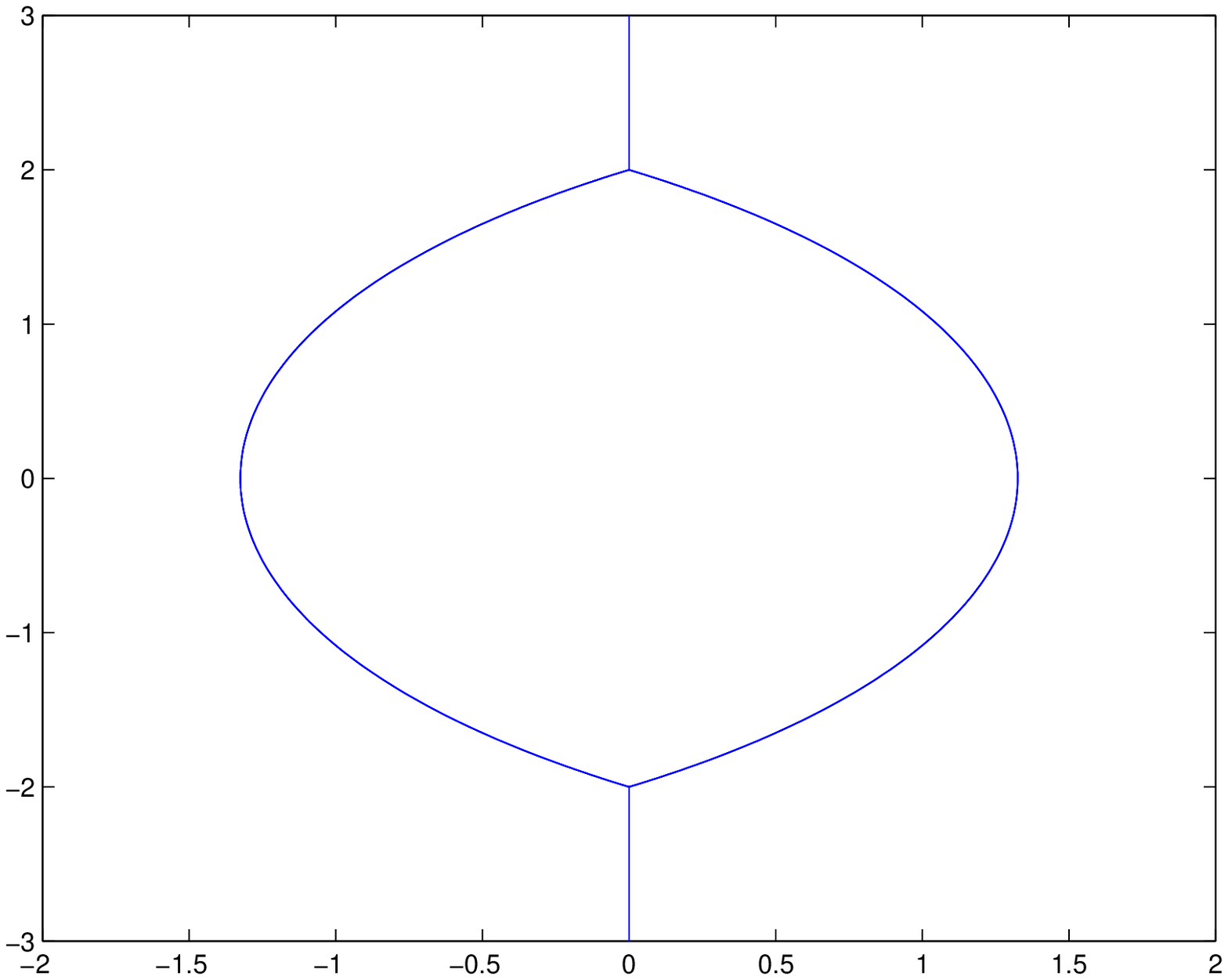}
}
\caption{Curves where $\frac{1}{2\pi i} \int_{\bar\beta}^z \frac{R(s)}{s}\, ds$
is real, for the values $A = 1.01$ (left) and $A = 2$ (right). The two points of intersection
are $\bar\beta$ and $\beta$.}
\label{fig:plot2}
\end{figure}

A major step in the analysis of the Riemann--Hilbert problem
for $Y$ is the selection of the right contour. In order that the
subsequent analysis works, the contour cannot be arbitrary but has to chosen
in a precise way. The contour depends on $A$.
We define
\begin{equation} \label{eq31}
    \beta = 2 - A + 2 i \sqrt{A-1}
\end{equation}
Martinez et al.\ \cite{MMO1} showed that the values of $z$ for which
\[ \frac{1}{2\pi i} \int_{\bar\beta}^z \frac{(s-\beta)^{1/2} (s-\bar\beta)^{1/2}}{s}\, ds
    \quad \mbox{ is real}, \]
where the branch of the square roots is chosen so that they are analytic and single valued
on the path of integration from $\bar\beta$ to $z$, form a system of curves as shown in
Figure \ref{fig:plot2} for a number of values of $A$.
In geometric function theory these curves are known as trajectories of the quadratic
differential $\frac{(s-\beta)(s-\bar\beta)}{s^2}\, d^2 s$, see \cite{Strebel}.
We see four smooth (in fact analytic) curves. Two curves are connecting $\bar\beta$
with $\beta$, one of them crosses the negative real axis, and the other one
crosses the positive real axis.
\begin{definition}
We define $\Gamma$ as the trajectory of the quadratic differential
$\frac{(s-\beta)(s-\bar\beta)}{s^2}\, d^2 s$ from $\bar\beta$ to $\beta$
which crosses the negative real axis. $\Gamma$ is oriented from $\bar\beta$
to $\beta$.
\end{definition}

We  put
\begin{equation} \label{eq32}
    R(z) = (z-\beta)^{1/2} (z-\bar{\beta})^{1/2}, \qquad z \in \mathbb C \setminus \Gamma,
\end{equation}
where the branch is chosen which is defined and analytic on $\mathbb C \setminus \Gamma$,
and which is such that $R(z) \sim z$ as $z \to \infty$. For $s \in \Gamma$, we use
$R_+(s)$ and $R_-(s)$ to denote the limits from the $+$-sides and
$-$-sides, respectively. As usual, the $+$-side of an oriented curve lies to the left, and
the $-$-side lies to the right, if one traverses the curve.

Then by definition of $\Gamma$, we have
\begin{equation} \label{eq33}
    \frac{1}{2\pi i} \int_{\bar\beta}^z \frac{R_+(s)}{s} \, ds
    \quad \mbox{ is real for every } z \in \Gamma,
\end{equation}
with integration along the $+$-side of $\Gamma$.

\begin{figure}[th!]
\centerline{\includegraphics[width=5cm]{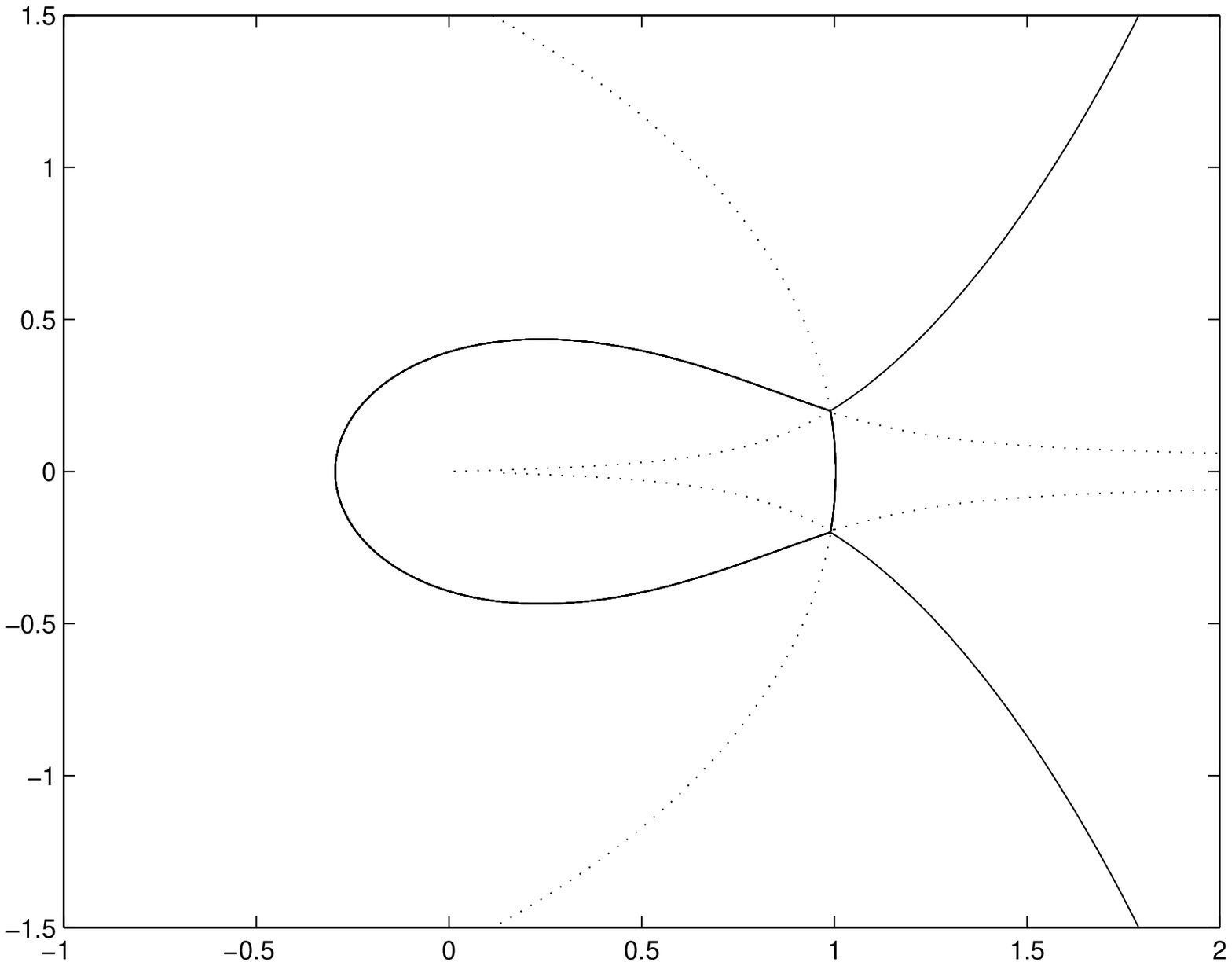} \qquad \qquad
    \includegraphics[width=5cm]{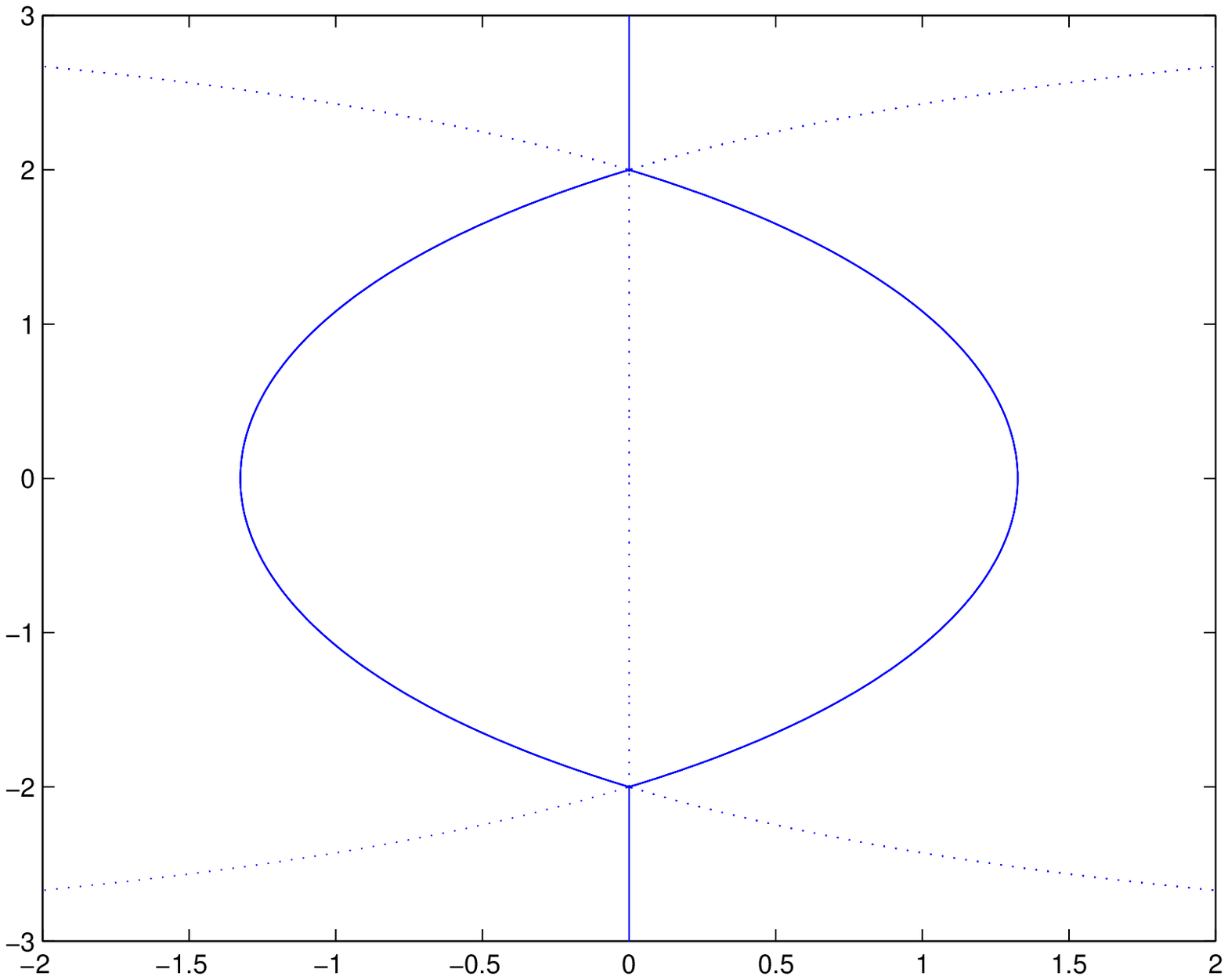}
}
\caption{Curves where $\frac{1}{2\pi i} \int_{\bar\beta}^z \frac{R(s)}{s}\, ds$ is real
(solid lines) and curves where $\frac{1}{2\pi i} \int_{\bar\beta}^z \frac{R(s)}{s}\, ds$
or $\frac{1}{2\pi i} \int_{\beta}^z \frac{R(s)}{s}\, ds$ is purely imaginary (dotted lines), for
the values $A = 1.01$ (left) and $A = 2$ (right).}
\label{fig:plot3}
\end{figure}

It is also of interest to know where
$\frac{1}{2 \pi i} \int_{\bar\beta}^z \frac{R(s)}{s}\, ds$
or $\frac{1}{2\pi i} \int_{\beta}^z \frac{R(s)}{s}\, ds$
is purely imaginary. These are the dotted curves shown in Figure \ref{fig:plot3}.
The dotted curves are analytic extensions of the solid ones.

We can now state which contour $\Sigma$ to choose.
\begin{definition}
We let $\Sigma$ be the contour in $\cal F$ consisting of $\Gamma$ together
with the two dotted curves that form the analytic extension of $\Gamma$.

We denote the part of $\Sigma \setminus \Gamma$ in the lower half-plane
by $\Sigma_1$ and its mirror image in the upper half-plane  by $\Sigma_2$.
\end{definition}

\begin{figure}[th!]
\centerline{\includegraphics[width=5cm]{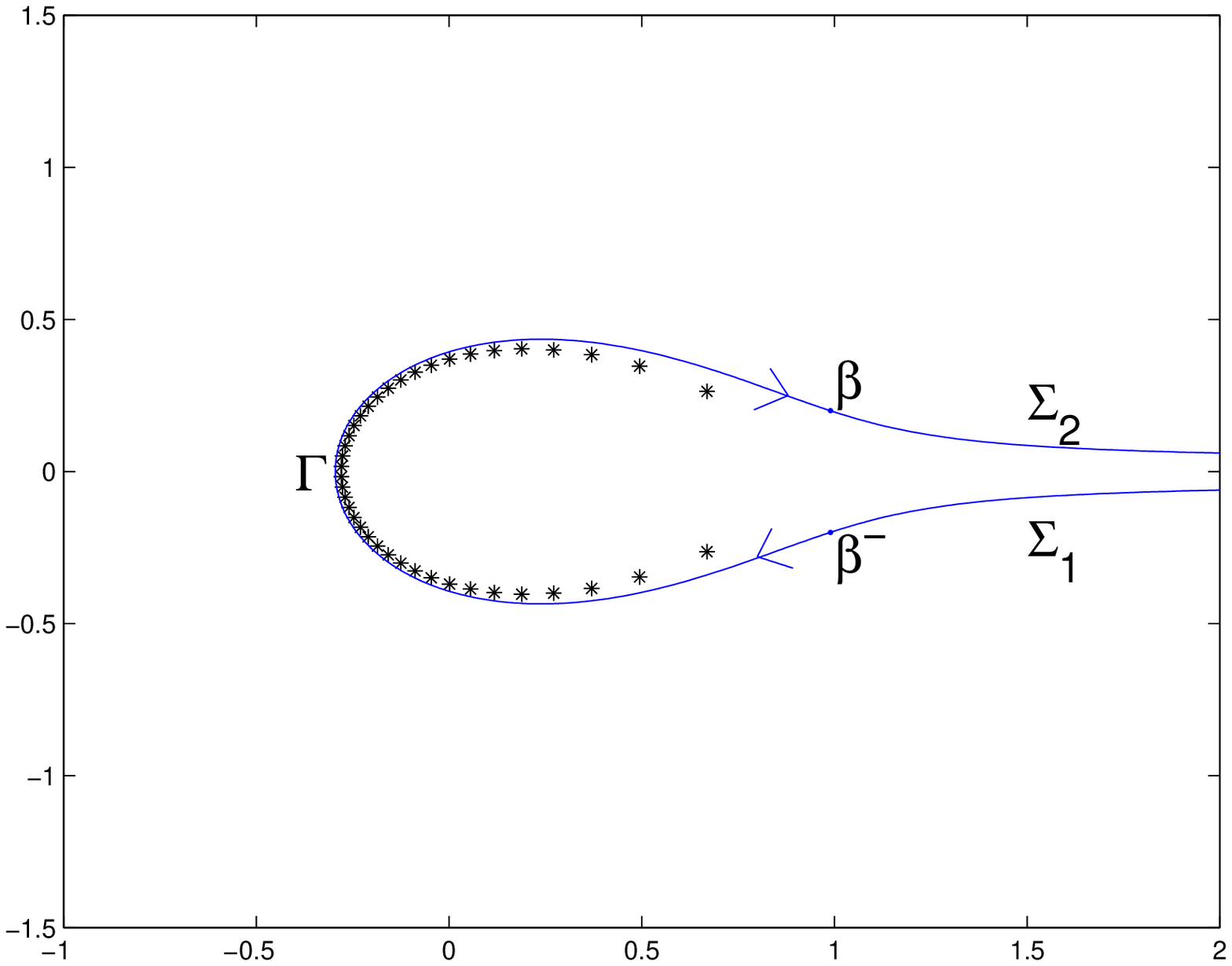} \qquad \qquad
    \includegraphics[width=5cm]{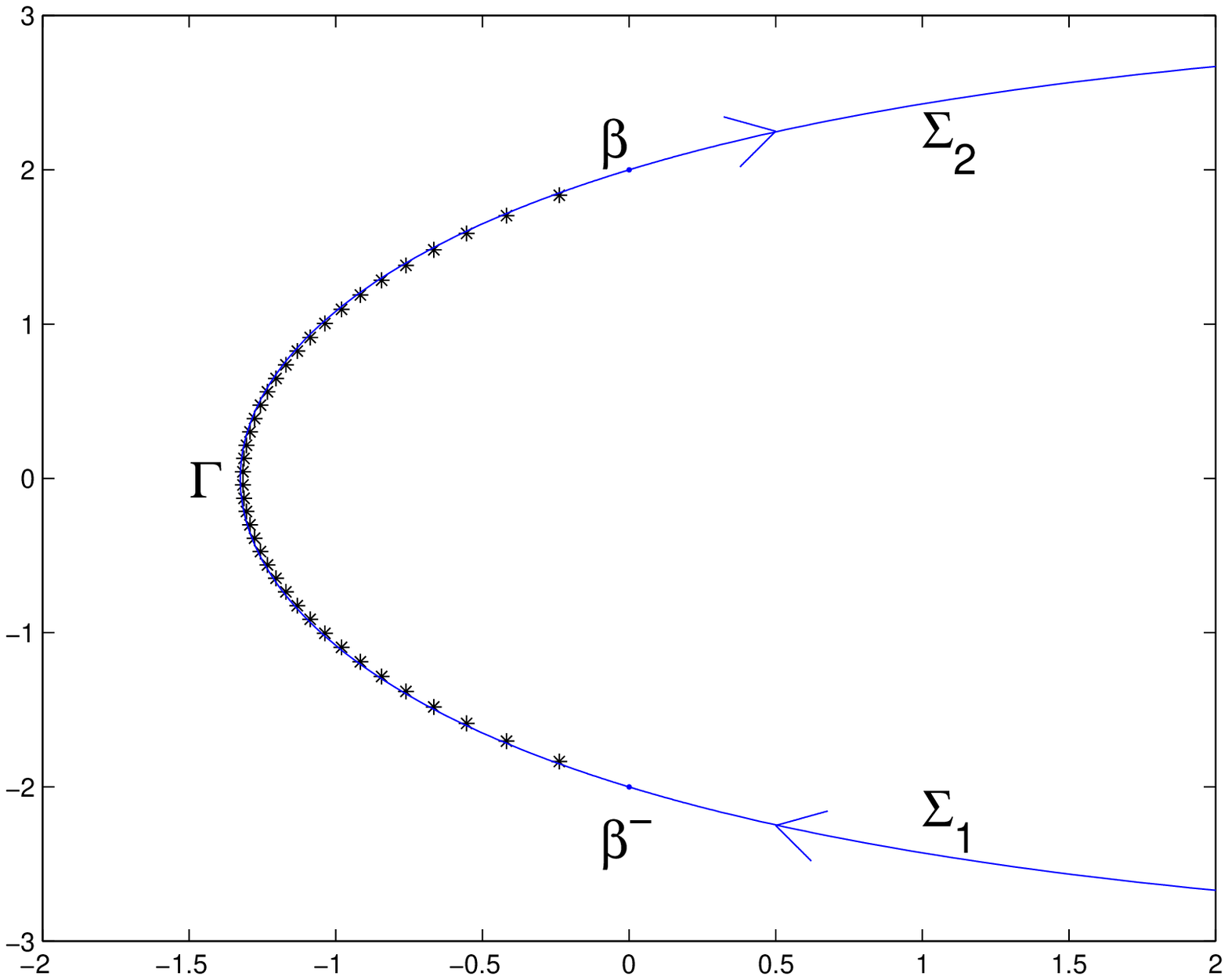}
}
\caption{The curves $\Sigma = \Gamma \cup \Sigma_1 \cup \Sigma_2$,
and the zeros of $L_{40}^{(-40A)}(40x)$, for the values $A = 1.01$ (left)
and $A = 2$ (right).}
\label{fig:plot4}
\end{figure}

So we have a disjoint union $\Sigma = \Gamma \cup \Sigma_1 \cup \Sigma_2$.
Figure \ref{fig:plot4} shows the curve $\Sigma$ for two values of $A$,
together with the zeros of the corresponding Laguerre polynomial
$L_{40}^{(-40A)}(40x)$ of degree $40$.
The figure shows that the zeros are close to $\Gamma$, and that they
are in the domain $\Omega_-$. These findings will be confirmed by our
final result, Corollary 7.2 below.

\begin{remark}
We have chosen $\Sigma_2$ so that $\int_{\beta}^z \frac{R(s)}{s}\, ds$
is real and positive on $\Sigma_2$. This is not essential.
What is important for the subsequent
analysis is that it has positive real part on $\Sigma_2$. This means that we have the
freedom to deform  $\Sigma_2$, as long as we take care that the real part
of $\int_{\beta}^z \frac{R(s)}{s}\, ds$ is positve on $\Sigma_2$.
However, it will be convenient to do this deformation only away from $\beta$,
so that in a neighborhood of $\beta$, we have $\Sigma_2$ exactly as we
defined it.

Similar remarks apply to $\int_{\bar\beta}^z \frac{R(s)}{s}\, ds$ and $\Sigma_1$.
\end{remark}

\subsection{A probability measure on $\Gamma$}
The following proposition gives one of the crucial properties of $\Gamma$.

\begin{proposition}
The (complex) measure $\frac{1}{2\pi i} \frac{R_+(s)}{s}\, ds$ is a probability
measure on $\Gamma$.
\end{proposition}
\begin{proof}
We show first that
\begin{equation} \label{eq34}
    \frac{1}{2\pi i} \int_{\Gamma} \frac{R_+(s)}{s}\, ds = 1.
\end{equation}
Let $I = \int_{\Gamma} \frac{R_+(s)}{s}\, ds$.
Since $R_-(s) = - R_+(s)$ for $s \in \Gamma$, we have
\[ 2 I = \int_{\Gamma} \left( \frac{R_+(s)}{s} - \frac{R_-(s)}{s} \right)\, ds
    = \oint_{\gamma} \frac{R(s)}{s}\, ds, \]
where $\gamma$ is a closed contour encircling the curve $\Gamma$ once in the clockwise
direction and not encircling $z=0$.

After contour deformation, we  pick up residues at $z=0$ and at $z=\infty$,
namely
\begin{equation} \label{eq35}
2I = 2\pi i \Res_{z=0} \frac{R(z)}{z} - 2\pi i \Res_{z =\infty} \frac{R(z)}{z}.
\end{equation}
The residue at $z=0$ is
\begin{equation} \label{eq36}
\Res_{z=0} \frac{R(z)}{z} = R(0) = |\beta| = \sqrt{(2-A)^2 + 4(A-1)} = A
\end{equation}
and the residue at $z = \infty$ is the coefficient of $z^{-1}$ in the
Laurent expansion of $R(z)/z$
\begin{eqnarray*}
    \frac{R(z)}{z} & = & (1-\beta/z)^{1/2} (1-\bar{\beta}/z)^{1/2} \\
    & = & \left(1- \frac{\beta}{2z} + O(z^{-2})\right)
        \left(1-\frac{\bar{\beta}}{2z} + O(z^{-2})\right) \\
    & = & 1 - \frac{\beta + \bar{\beta}}{2z} + O(z^{-2}).
\end{eqnarray*}
Thus
\begin{equation} \label{eq37}
    \Res_{z=\infty} \frac{R(z)}{z} = - \frac{\beta + \bar{\beta}}{2} =
    - \Re \beta = - (2-A).
\end{equation}
Hence by (\ref{eq35})--(\ref{eq37}) we have $2 I = 2\pi i (A  + (2-A)) = 4 \pi i$,
so that (\ref{eq34}) follows.

Having (\ref{eq34}) we can now prove the proposition. Let $t \mapsto z(t)$
for $t \in [0,t_0]$, be the arc length parametrization of $\Gamma$
starting at $\bar\beta$. Thus $z(0) = \bar\beta$ and $z(t_0) = \beta$.
Then we have for $z = z(t)$ with $t \in (0,t_0)$,
\begin{equation} \label{eq38}
    \frac{1}{2\pi i} \int_{\bar\beta}^z \frac{R_+(s)}{s}\, ds =
    \frac{1}{2\pi i} \int_0^t \frac{R_+(z(\tau))}{z(\tau)} z'(\tau)\, d\tau,
\end{equation}
and this is real for every $t$ by construction of $\Gamma$.
It has the value $0$ for $t = 0$ and the value $1$ for $t = t_0$
(due to (\ref{eq34})).
The derivative of (\ref{eq38}) with respect to $t$ is
$\frac{1}{2\pi i} \frac{R_+(z(t))}{z(t)} z'(t)$
which is not zero for $t \in (0,t_0)$. Thus (\ref{eq38}) can only increase
from $0$ to $1$ as $t$ increases from $0$ to $t_0$.
Hence  $\frac{1}{2\pi i}  \frac{R_+(s)}{s}\, ds$ is a positive measure
on $\Gamma$. It is a probability measure because of (\ref{eq34}).
\end{proof}

\subsection{Auxiliary functions}

With the measure $d\mu(s) = \frac{1}{2\pi i} \frac{R_+(s)}{s} ds $ on $\Gamma$
we define the so-called $g$-function as follows.
\begin{definition}
The $g$-function is the complex logarithmic potential of $\mu$, that is,
\begin{equation} \label{eq39}
    g(z) = \int_{\Gamma} \log(z-s) \, d\mu(s),
      \qquad z \in \mathbb C \setminus (\Gamma \cup \Sigma_1),
\end{equation}
where for each $s$ we view $\log (z-s)$ as an analytic function of the variable
$z$, with branch cut emanating from $z=s$. The cut is taken
along $\Gamma \cup \Sigma_1$.
\end{definition}

We need two more functions.
\begin{definition}
The $\phi$-function is defined as
\begin{equation} \label{eq310}
\phi(z) = \frac{1}{2} \int_{\beta}^z \frac{R(s)}{s} \, ds,
    \qquad z \in \mathbb C \setminus (\Gamma \cup \Sigma_1 \cup [0, \infty)),
\end{equation}
where the path of integration from $\beta$ to $z$ lies entirely
in $\mathbb C \setminus (\Gamma \cup \Sigma_1 \cup [0,\infty))$,
except for the initial point $\beta$.

The $\tilde{\phi}$-function is defined as
\begin{equation} \label{eq311}
\tilde{\phi}(z) = \frac{1}{2} \int_{\bar\beta}^z \frac{R(s)}{s} \, ds,
    \qquad z \in \mathbb C \setminus (\Gamma \cup \Sigma_2 \cup [0, \infty)),
\end{equation}
where the path of integration from $\bar\beta$ to $z$ lies entirely
in $\mathbb C \setminus (\Gamma \cup \Sigma_2 \cup [0,\infty))$,
except for the initial point $\bar\beta$.
\end{definition}

It is immediate from (\ref{eq33}) that
$\tilde{\phi}_+(z)$ is purely imaginary for $z \in \Gamma$.
By Proposition 3.4 its imaginary part increases from $0$ to $\pi$ as $z$ traverses
the curve $\Gamma$ from $\bar\beta$ to $\beta$.
In particular we have $\tilde{\phi}_+(\beta) = \pi i$. Similarly
$\tilde{\phi}_-(\beta) = - \pi i$. Thus we have
\begin{equation} \label{eq312}
    \tilde{\phi}(z) =  \left\{ \begin{array}{ll}
        \phi(z) + \pi i & \quad \mbox{for } z \in \Omega_+, \\[10pt]
        \phi(z) - \pi i & \quad \mbox{for } z \in \Omega_-.
    \end{array} \right.
\end{equation}

\begin{proposition}
There is a constant $\ell$ such that
\begin{equation} \label{eq313}
    g(z) = \frac{1}{2} \left( A \log z + z + \ell \right) - \phi(z),
    \qquad z \in \mathbb C \setminus (\Gamma \cup \Sigma_1 \cup [0,\infty)),
\end{equation}
where $\log z$ is defined with a branch cut along $[0,\infty)$.
\end{proposition}
\begin{proof}
We note that
\[ g'(z) = \frac{1}{2\pi i} \int_{\Gamma} \frac{1}{z-s} \frac{R_+(s)}{s}\, ds
    = \frac{1}{2} \frac{1}{2\pi i} \oint_{\gamma} \frac{1}{z-s} \frac{R(s)}{s}\, ds \]
where $\gamma$ is a closed contour in $\mathbb C \setminus \Gamma$, that
encircles $\Gamma$ once in the clockwise direction, but does not encircle $z$ and $0$.
Then exactly as in the proof of Proposition 3.4, we deform the contour, and now pick
up residues at $z$, $0$ and $\infty$. We find for $z \in \mathbb C \setminus \Gamma$,
\begin{eqnarray} \nonumber
g'(z) & = & \frac{1}{2} \left[ \Res_{s=z} \left(  \frac{1}{z-s} \frac{R(s)}{s} \right) +
    \Res_{s=0} \left( \frac{1}{z-s} \frac{R(s)}{s} \right) -
    \Res_{s=\infty} \left( \frac{1}{z-s} \frac{R(s)}{s} \right)
    \right] \\ \nonumber
    & = & \frac{1}{2} \left[ - \frac{R(z)}{z} + \frac{R(0)}{z} + 1 \right] \\
    \label{eq314}
    & = & \frac{1}{2} \left[ - \frac{R(z)}{z} + \frac{A}{z} + 1 \right].
\end{eqnarray}

Let $z \in \mathbb C \setminus (\Gamma \cup \Sigma_1 \cup [0,\infty))$.
Integrating (\ref{eq314}) from $\beta$ to $z$ along a curve in
$\mathbb C \setminus (\Gamma \cup \Sigma_1 \cup [0,\infty)$, we find
\begin{equation} \label{eq315}
    g(z) = g(\beta) + \frac{1}{2} (A \log z + z) - \frac{1}{2} (A \log \beta + \beta)
    - \phi(z).
\end{equation}
Thus (\ref{eq313}) holds with $\ell = 2 g(\beta) - (A \log \beta + \beta)$.
\end{proof}

\subsection{Jump properties of $g$}
From Proposition 3.7 we obtain the following jump relations for $g$ across
the contour $\Sigma$. These jumps are crucial for the subsequent analysis.

\begin{proposition}
\begin{enumerate}
\item[\rm (a)] We have
\begin{equation} \label{eq316}
    g_+(z) - g_-(z) = 2 \pi i  \quad \mbox{for } z \in \Sigma_1,
\end{equation}
and
\begin{equation} \label{eq317}
    g_+(z) - g_-(z) = - \phi_+(z) + \phi_-(z) = -2\phi_+(z) = 2 \phi_-(z)
    \qquad \mbox{for } z \in \Gamma.
\end{equation}
\item[\rm (b)] We have, with the same constant $\ell$ as in Proposition {\rm 3.7},
\begin{equation} \label{eq318}
     g_+(z) + g_-(z) = A \log z + z + \ell \qquad \mbox{ for } z \in \Gamma,
\end{equation}
\begin{equation} \label{eq319}
    g_+(z) + g_-(z) = A \log z + z + \ell - 2\tilde{\phi}(z)
    \qquad \mbox{ for } z \in \Sigma_1,
\end{equation}
and
\begin{equation} \label{eq320}
    g_+(z) + g_-(z) = A \log z + z + \ell - 2\phi(z)
    \qquad \mbox{ for } z \in \Sigma_2.
\end{equation}
\end{enumerate}
\end{proposition}
\begin{proof}
In (\ref{eq313}) we let $z$ approach $\Sigma$, from the $+$- and $-$-sides,
respectively, to obtain
\begin{equation} \label{eq321}
    g_{\pm}(z) = \frac{1}{2} (A \log z + z) +  \frac{1}{2} \ell -
    \phi_{\pm}(z),
    \qquad \mbox{for } z \in \Sigma.
\end{equation}
Since $\phi$ changes sign across $\Gamma$, (\ref{eq317}) and (\ref{eq318})
immediately follow from (\ref{eq321}).
For $z \in \Sigma_1$, we have by (\ref{eq312}) and (\ref{eq321})
\begin{eqnarray*}
    g_+(z) - g_-(z) & = & -\phi_+(z) + \phi_-(z) \ = \
        -(\tilde{\phi}_+(z) - \pi i) + (\tilde{\phi}_-(z) + \pi i) \\
        & = & 2 \pi i - \tilde{\phi}_+(z) + \tilde{\phi}_-(z),
\end{eqnarray*}
which is (\ref{eq316}) as $\tilde{\phi}$ is analytic across $\Sigma_1$.
In addition, we have
\begin{eqnarray*}
    g_+(z) + g_-(z) & = & A \log z + z + \ell - \phi_+(z) - \phi_-(z) \\[10pt]
    & = & A \log z + z + \ell - \tilde{\phi}_+(z) - \tilde{\phi}_-(z)
    \qquad \mbox{ for } z \in \Sigma_1.
\end{eqnarray*}
which yields (\ref{eq319}). Similarly, (\ref{eq320}) follows.
\end{proof}

\section{First two transformations: $Y \mapsto U \mapsto T$}
\setcounter{equation}{0}

\subsection{First transformation $Y \mapsto U$}
With the $g$-function and the constant $\ell$ from Proposition 3.7,
we perform the first transformation of the Riemann--Hilbert problem.

\begin{definition}
We define for
$z \in \mathbb C \setminus \Sigma$,
\begin{equation} \label{eq41}
    U(z) = e^{-n (\ell/2) \sigma_3} Y(z)
    e^{-ng(z) \sigma_3} e^{n (\ell/2) \sigma_3}.
\end{equation}
Here, and in what follows, $\sigma_3$ denotes the Pauli matrix
$\sigma_3 =\left( \begin{array}{cc} 1 & 0 \\ 0 & -1 \end{array} \right)$,
so that for example $e^{-ng(z) \sigma_3} = \left( \begin{array}{cc}
    e^{-ng(z)} & 0 \\ 0 & e^{ng(z)}
    \end{array} \right)$.
\end{definition}

From the Riemann--Hilbert problem for $Y$ it follows by a straightforward
calculation that $U$ is the unique solution of the following Riemann--Hilbert problem.

\subsubsection*{Riemann--Hilbert problem for $U$:}
The problem is to determine a $2\times 2$ matrix valued function
$U : \mathbb C \setminus \Sigma \to \mathbb C^{2\times 2}$ such that
\begin{enumerate}
\item[(a)] $U(z)$ is analytic for $z \in \mathbb C \setminus \Sigma$,
\item[(b)] $U(z)$ possesses continuous boundary values for $z \in \Sigma$,
denoted by $U_+(z)$ and $U_-(z)$,  and
\begin{equation} \label{eq42}
    U_+(z) = U_-(z)
    \left( \begin{array}{cc}
    e^{-n(g_+(z) -g_-(z))} & z^{-An} e^{-nz} e^{n(g_+(z)+g_-(z) - \ell)} \\
    0 & e^{n(g_+(z) - g_-(z))} \end{array} \right)
\end{equation}
for $z \in \Sigma$,
\item[(c)] $U(z)$ behaves like the identity at infinity:
\begin{equation} \label{eq43}
    U(z) = I + O\left(\frac{1}{z}\right)
    \qquad \mbox{ as } z \to \infty, \quad z \in \mathbb C \setminus \Sigma.
\end{equation}
\end{enumerate}

The jump relation (\ref{eq42}) for $U$ has a different form on
the three parts $\Gamma$, $\Sigma_1$ and $\Sigma_2$.
On $\Gamma$ we have that $g_+ - g_- = -2\phi_+ = 2\phi_-$ by (\ref{eq317}),
and $g_+ + g_- = A \log z + z + \ell$ by (\ref{eq318}) so that
\begin{equation} \label{eq44}
    U_+(z) = U_-(z)
    \left(\begin{array}{cc}
    e^{2n\phi_+(z)} & 1 \\ 0 & e^{2n \phi_-(z)}
    \end{array} \right)
    \qquad \mbox{for } z \in \Gamma.
\end{equation}
On $\Sigma_1$ we use (\ref{eq316}) and (\ref{eq319}) to obtain
\begin{equation} \label{eq45}
    U_+(z) = U_-(z)  \left( \begin{array}{cc}
    1 & e^{-2n \tilde{\phi}(z)}  \\
    0 & 1 \end{array} \right)
    \qquad \mbox{ for } z \in \Sigma_1,
\end{equation}
and similarly it follows that
\begin{equation} \label{eq46}
    U_+(z) = U_-(z)  \left( \begin{array}{cc}
    1 & e^{-2n \phi(z)}  \\
    0 & 1 \end{array} \right)
    \qquad \mbox{ for } z \in \Sigma_2.
\end{equation}

The transformation $Y\mapsto U$ has the effect of normalizing the Riemann--Hilbert
problem at infinity. In addition, by the construction of $\Sigma$,
we have that $\phi$ is real and positive on $\Sigma_2$, and $\tilde{\phi}$ is
real and positive on $\Sigma_1$. So the jump matrices for $U$
in (\ref{eq45}) and (\ref{eq46}) are close to the identity if $n$ is large.
Since $\phi$ has purely imaginary boundary values on both sides of $\Gamma$,
the jump matrix for $U$ on $\Gamma$ in (\ref{eq44})
has oscillatory diagonal entries.

\subsection{Second transformation $U \mapsto T$}

\begin{figure}[th!]
\centerline{\includegraphics*[width=10.8cm,height=6.4cm]{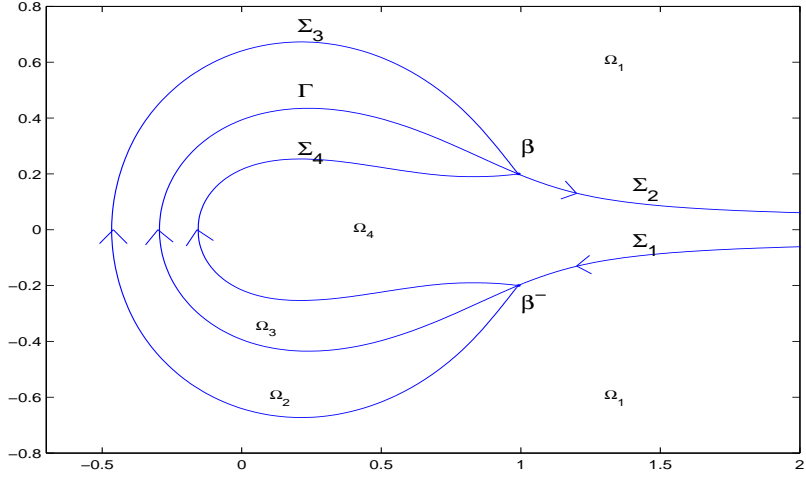}}
\caption{Contour $\Sigma^T = \Gamma \cup \bigcup_j \Sigma_j$ and the domains $\Omega_j$,
$j=1,\ldots, 4$, for the Riemann--Hilbert problem for $T$, for the value $A = 1.01$.}
\label{fig:plotT}
\end{figure}

The jump matrix for $U$ on $\Gamma$, see (\ref{eq44}), factors as
\begin{equation} \label{eq47}
    \left(\begin{array}{cc}
    e^{2n\phi_+(z)} & 1 \\ 0 & e^{2n \phi_-(z)}
    \end{array} \right) =
    \left(\begin{array}{cc}
    1 & 0 \\ e^{2n \phi_-(z)} & 1
    \end{array} \right)
    \left(\begin{array}{cc}
    0 & 1 \\ -1 & 0
    \end{array} \right)
    \left(\begin{array}{cc}
    1 & 0 \\ e^{2n \phi_+(z)} & 1
    \end{array} \right).
\end{equation}
Observe that the first matrix in the right-hand side of (\ref{eq47})
can be analytically continued to the $-$-side of the contour $\Gamma$,
and in doing so, the (1,2) entry becomes exponentially decaying in $n$.
Similarly, the third matrix in the right-hand side of (\ref{eq47}) can be
analytically continued to the $+$-side of the contour $\Gamma$,
and in doing so, the (1,2) entry also becomes exponentially decaying in $n$.
We are thus led to introduce the following ``contour augmentation'' step,
as part of the steepest descent / stationary phase method for Riemann--Hilbert
problems developed by Deift and Zhou.
The oriented contour $\Sigma^T$ consists of $\Sigma$ plus two simple curves
$\Sigma_3$ and $\Sigma_4$ from $\bar\beta$ to $\beta$, contained in $\Omega_+$ and $\Omega_-$,
respectively, as shown in Figure \ref{fig:plotT}.
We choose $\Sigma_3$ and $\Sigma_4$ such that $\Re \phi(z) < 0$ on $\Sigma_3$ and $\Sigma_4$.

It is possible to choose such curves. Indeed, $\phi$ is positive on $\Sigma_2$, and
its real part vanishes on $\Gamma$ and on the other solid lines shown in Figure
\ref{fig:plot3}.
So $\Re \phi > 0$ in the full region on the right. In the two other regions, bounded
by the solid lines, we then have that $\Re \phi < 0$. Note that $\Re \phi$ does not
change sign across $\Gamma$.

Note that by (\ref{eq312}) we also have $\Re \tilde{\phi}(z) < 0$ on
$\Sigma_3$ and $\Sigma_4$.

Then $\mathbb C \setminus \Sigma^T$ has four connected components, denoted by
$\Omega_1$, $\Omega_2$, $\Omega_3$, and $\Omega_4$ as indicated in Figure \ref{fig:plotT}.

\begin{definition}
We define $T : \mathbb C \setminus \Sigma^T \to \mathbb C^{2 \times 2}$ by
\begin{equation} \label{eq48}
    T(z) = U(z) \qquad \mbox{ for } z \in \Omega_1 \cup \Omega_4,
\end{equation}
\begin{equation} \label{eq49}
    T(z) = U(z)  \left(\begin{array}{cc} 1 & 0 \\
    -e^{2n \phi(z)} & 1 \end{array} \right)
    \qquad \mbox{ for } z \in \Omega_2,
\end{equation}
\begin{equation} \label{eq410}
    T(z) = U(z)  \left(\begin{array}{cc} 1 & 0 \\
    e^{2n \phi(z)} & 1 \end{array} \right)
    \qquad \mbox{ for } z \in \Omega_3,
\end{equation}
\end{definition}

Then from the Riemann--Hilbert problem for $U$ and the factorization (\ref{eq47})
we obtain that $T$ is the unique solution of the following Riemann--Hilbert problem.

\subsubsection*{Riemann--Hilbert problem for $T$:}
The problem is to determine a $2\times 2$ matrix valued function
$T : \mathbb C \setminus \Sigma^T \to \mathbb C^{2\times 2}$ such that
the following hold:
\begin{enumerate}
\item[(a)] $T(z)$ is analytic for $z \in \mathbb C \setminus \Sigma^T$,
\item[(b)] $T(z)$ possesses continuous boundary values for $z \in \Sigma^T$,
denoted by $T_+(z)$ and $T_-(z)$, and
\begin{equation} \label{eq411}
    T_+(z) = T_-(z)
   \left(\begin{array}{cc} 0 & 1 \\ -1 & 0 \end{array} \right)
   \qquad \mbox{ for } z\in \Gamma,
\end{equation}
\begin{equation} \label{eq412}
    T_+(z) = T_-(z) \left(\begin{array}{cc} 1 & 0 \\
    e^{2n \phi(z)} & 1 \end{array} \right)
    \qquad \mbox{ for } z \in \Sigma_3 \cup \Sigma_4,
\end{equation}
\begin{equation} \label{eq413}
    T_+(z) = T_-(z) \left( \begin{array}{cc}
     1 & e^{-2n \tilde{\phi}(z)} \\
     0 & 1 \end{array} \right)
     \qquad \mbox{ for } z\in \Sigma_1,
\end{equation}
and
\begin{equation} \label{eq414}
    T_+(z) = T_-(z) \left( \begin{array}{cc}
     1 & e^{-2n \phi(z)} \\
     0 & 1 \end{array} \right)
     \qquad \mbox{ for } z\in \Sigma_2,
\end{equation}
\item[(c)] $T(z)$ behaves like the identity at infinity:
\begin{equation} \label{eq415}
    T(z) = I + O\left(\frac{1}{z}\right)
    \qquad \mbox{ as } z \to \infty, \quad z \in \mathbb C \setminus \Sigma^T.
\end{equation}
\end{enumerate}

\section{Construction of the parametrix for $T$}
\setcounter{equation}{0}

\subsection{Parametrix away from $\beta$ and $\bar\beta$}
As remarked following (\ref{eq47}), the jump matrices  appearing in
(\ref{eq412}) for $z \in \Sigma_3 \cup \Sigma_4$ are exponentially close
to the identity matrix away from $\beta$ and $\bar\beta$. Similarly,
the jump matrices appearing in (\ref{eq413}) and (\ref{eq414})
are  also exponentially close to the identity matrix away from $\beta$ and
$\bar\beta$. This hints that these portions of the contour on which the
Riemann--Hilbert problem for $T$ is posed should be somehow negligible.
Thus we expect that the leading order asymptotics is determined by the
solution $N$
of the following model Riemann--Hilbert problem:

\subsubsection*{Riemann--Hilbert problem for $N$:}
The problem is to determine $N : \mathbb C \setminus \Gamma \to \mathbb C^{2\times 2}$
such that the following hold.
\begin{enumerate}
\item[(a)] $N(z)$ is analytic for $z \in \mathbb C \setminus \Gamma$,
\item[(b)] $N(z)$ possesses continuous boundary values for $z \in
\Gamma \setminus \{\beta, \bar\beta\}$,
denoted by $N_+(z)$ and $N_-(z)$, and
\begin{equation} \label{eq51}
    N_+(z) = N_-(z)
    \left(\begin{array}{cc} 0 & 1 \\ -1 & 0 \end{array} \right)
    \qquad \mbox{ for } z \in \Gamma \setminus \{\beta, \bar\beta\},
\end{equation}
\item[(c)] $N(z) = I + O\left(\frac{1}{z}\right)$ for $z \to \infty$.
\end{enumerate}

This Riemann--Hilbert problem for $N$ is solved explicitly
by, see \cite[p.1520]{DKMVZ2} or \cite[p.200]{Deift},
\begin{equation} \label{eq52}
    N(z) =
    \left(\begin{array}{cc}
    \frac{a(z) + a(z)^{-1}}{2} &
    \frac{a(z) - a(z)^{-1}}{2i} \\[0.5pt]
    \frac{a(z) - a(z)^{-1}}{-2i} &
    \frac{a(z) + a(z)^{-1}}{2}
    \end{array} \right),
\end{equation}
where
\begin{equation} \label{eq53}
    a(z) =  \frac{(z-\beta)^{1/4}}{(z-\bar\beta)^{1/4}}.
\end{equation}
The branches of the roots in (\ref{eq53}) are chosen such that $a(z)$ is analytic
on $\mathbb C \setminus \Gamma$ and $\lim\limits_{z \to \infty} a(z) = 1$.

\begin{remark}
The solution (\ref{eq52}) is not the only solution to the Riemann--Hilbert problem
for $N$. It is the unique solution that satisfies, in addition to (a), (b), and (c),
the condition
\begin{enumerate}
\item[(d)] Near the endpoints $\beta$ and $\bar\beta$, we have
\[ N(z) = O(|z-\beta|^{-1/4}) \qquad \mbox{ as } z \to \beta, \]
\[ N(z) = O(|z-\bar{\beta}|^{-1/4}) \qquad \mbox{ as } z \to \bar\beta, \]
with the $O$-term being taken entry-wise.
\end{enumerate}
\end{remark}

\begin{remark}
For explicit calculations, it is useful to have an alternative
expression for $N$. For $z = x$ real, it is clear from (\ref{eq53})
that $a(x)$ has modulus one.  If $\arg(a(x)) = \theta(x)$, then (\ref{eq52})
shows
\[ N(x) = \left(\begin{array}{rc}
    \cos \theta(x) & \sin \theta(x) \\
    -\sin \theta(x) & \cos \theta(x)
    \end{array} \right). \]
Since $\tan ( 2\theta(x)) =  -\frac{2 \sqrt{A-1}}{x-2+A}$, we then find
\begin{equation} \label{eq54}
    N(z) = \left(\begin{array}{rc}
     \cos( \frac{1}{2} \arctan (\frac{2 \sqrt{A-1}}{z-2+A})) &
     - \sin( \frac{1}{2} \arctan (\frac{2 \sqrt{A-1}}{z-2+A})) \\[10pt]
     \sin( \frac{1}{2} \arctan (\frac{2 \sqrt{A-1}}{z-2+A})) &
     \cos( \frac{1}{2} \arctan (\frac{2 \sqrt{A-1}}{z-2+A}))
     \end{array}\right)
\end{equation}
first for $z=x$ real, but then also for arbitrary $z \in \mathbb C \setminus \Gamma$
by analytic continuation.
We have to take the appropriate branch of the multivalued $\arctan$ function in (\ref{eq54}).
Using trigonometric identities, one may then check that (\ref{eq54}) reduces to
\begin{equation} \label{eq55}
    N(z) = \left(\begin{array}{rc}
        \left( \frac{1 + R'(z)}{2} \right)^{1/2} &
        -\left(\frac{1-R'(z)}{2} \right)^{1/2} \\[10pt]
        \left(\frac{1-R'(z)}{2} \right)^{1/2} &
     \left( \frac{1 + R'(z)}{2} \right)^{1/2}
     \end{array} \right)
     \qquad \mbox{ for } z \in \mathbb C \setminus \Gamma,
\end{equation}
where as usual we have $R(z) = (z-\beta)^{1/2}(z-\bar\beta)^{1/2}$.
We can also verify directly that (\ref{eq55}) solves the Riemann--Hilbert problem
for $N$.
\end{remark}

\subsection{Parametrix near $\beta$}

\begin{figure}[th!]
\centerline{\includegraphics[width=10cm]{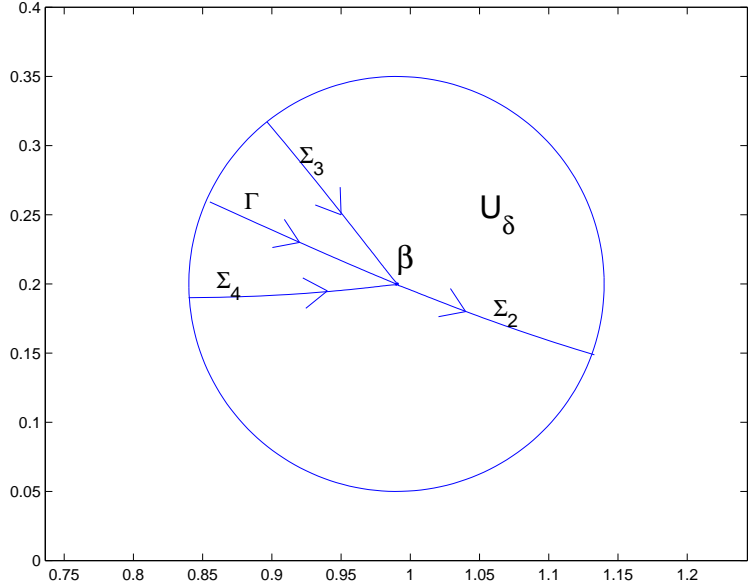}}
\caption{Neighborhood $U_{\delta}$ of $\beta$ and the parts of the contours $\Gamma$,
$\Sigma_2$, $\Sigma_3$, and $\Sigma_4$ that are within $U_{\delta}$.
The value of $A$ is $1.01$.}
\label{fig:plotP}
\end{figure}

The next step is a local analysis around the points $\beta$ and $\bar\beta$.
We need to construct a local parametrix $P$ in a neighborhood
$U_{\delta} = \{ z \in \mathbb C \mid |z-\beta| < \delta\}$ of $\beta$ such that
\begin{itemize}
\item $P$ satisfies the jumps for $T$ exactly in $U_{\delta}$,
\item $P$ matches $N$ on the boundary of $U_{\delta}$ up to order $1/n$.
\end{itemize}
See Figure \ref{fig:plotP} for the contours $\Sigma^T \cap U_{\delta}$.

More precisely, we have
\subsubsection*{Riemann--Hilbert problem for $P$:}
The problem is to determine, for a given $\delta > 0$ sufficiently small,
a $2\times 2$ matrix valued function
$P : \overline{U}_{\delta} \setminus \Sigma^T \to \mathbb C^{2\times 2}$
such that
\begin{enumerate}
\item[(a)] $P(z)$ is analytic for $z \in U_{\delta} \setminus \Sigma^T$,
and continuous on $\overline{U}_{\delta} \setminus \Sigma^T$,
\item[(b)] $P(z)$ possesses continuous boundary values for $z \in  \Sigma^T \cap U_{\delta}$,
denoted by $P_+(z)$ and $P_-(z)$, and
\begin{equation} \label{eq56}
    P_+(z) = P_-(z)
   \left(\begin{array}{cc} 0 & 1 \\ -1 & 0 \end{array} \right)
   \qquad \mbox{ for } z\in \Gamma \cap U_{\delta},
\end{equation}
\begin{equation} \label{eq57}
    P_+(z) = P_-(z) \left(\begin{array}{cc} 1 & 0 \\
    e^{2n \phi(z)} & 1 \end{array} \right)
    \qquad \mbox{ for } z \in (\Sigma_3 \cup \Sigma_4) \cap U_{\delta},
\end{equation}
and
\begin{equation} \label{eq58}
    P_+(z) = P_-(z) \left( \begin{array}{cc}
     1 & e^{-2n \phi(z)} \\
     0 & 1 \end{array} \right)
     \qquad \mbox{ for } z\in \Sigma_2 \cap U_{\delta},
\end{equation}
\item[(c)] There exists a constant $C > 0$ such that for every
$z \in \partial U_{\delta} \setminus \Sigma^T$,
\begin{equation} \label{eq59}
    \| P(z) N^{-1}(z)  - I \| \leq \frac{C}{n},
\end{equation}
    where $\| \cdot \|$ is any matrix norm.
\end{enumerate}

The construction of $P$ follows along the same lines as given by Deift et al.\ \cite{DKMVZ2}.
From its definition (\ref{eq310}) it is easy to see that the $\phi$-function has
a convergent expansion
\begin{equation} \label{eq510}
    \phi(z) = (z-\beta)^{3/2} \sum_{k=0}^{\infty} c_k (z-\beta)^k,
    \qquad c_0 \neq 0,
\end{equation}
in a neighborhood of $\beta$. The factor $(z-\beta)^{3/2}$ is defined with a cut
along $\Gamma \cup \Sigma_1$. Then
\begin{equation} \label{eq511}
    f(z) = \left[\frac{3}{2} \phi(z) \right]^{2/3}
\end{equation}
is defined and analytic in a neighborhood of $\beta$. We choose the $2/3$-root with
a cut along $\Gamma$ and such
that $f(z) > 0$ for $z \in \Sigma_2$. Recall $\phi > 0$ on $\Sigma_2$.

Then $f(\beta) = 0$ and $f'(\beta) \neq 0$. Therefore we can and do choose
$\delta$ so small that $\zeta = f(z)$ is a one-to-one mapping from $U_{\delta}$
onto a convex neighborhood $f(U_{\delta})$ of $\zeta = 0$.
Under the mapping $\zeta = f(z)$, we then have that
$\Sigma_2 \cap U_{\delta}$ corresponds to $(0,\infty) \cap f(U_{\delta})$
and that $\Gamma \cap U_{\delta}$ corresponds to $(-\infty,0] \cap f(U_{\delta})$.

Now we specify how to choose $\Sigma^T$ near $\beta$.
For an arbitrary, but fixed $\sigma \in (\pi/3, \pi)$, we choose $\Sigma_3$
and $\Sigma_4$ such that $f$ maps $\Sigma_3 \cap U_{\delta}$ and
$\Sigma_4 \cap U_{\delta}$ onto
$\{ \zeta \in f(U_{\delta}) \mid \arg \zeta = \sigma \}$
and
$\{ \zeta \in f(U_{\delta}) \mid \arg \zeta = - \sigma \}$,
respectively.

\begin{proposition}
The Riemann--Hilbert problem for $P$ is solved by
\begin{equation} \label{eq512}
    P(z) = E(z) \Psi^{\sigma}(n^{2/3} f(z)) e^{n\phi(z)\sigma_3}
\end{equation}
where
\begin{equation} \label{eq513}
    E(z) = \sqrt{\pi} e^{\frac{\pi i}{6}} \left(\begin{array}{cc}
        1 & -1 \\ -i & -i \end{array} \right)
        \left(\frac{n^{1/6} f(z)^{1/4}}{a(z)} \right)^{\sigma_3}
\end{equation}
and $\Psi^{\sigma}$ is an explicit matrix valued function built out of
the Airy function $\Ai$ and its derivative $\Ai'$ as follows
\begin{equation} \label{eq514}
    \Psi^{\sigma}(\zeta) = \left\{ \begin{array}{ll}
    \left(\begin{array}{cc}
    \Ai(\zeta) & \Ai(\omega^2 \zeta) \\
    \Ai'(\zeta) & \omega^2 \Ai'(\omega^2 \zeta) \end{array} \right)
    e^{- \frac{\pi i}{6} \sigma_3} &
    \mbox{for } 0 < \arg \zeta < \sigma, \\[10pt]
    \left(\begin{array}{cc}
    \Ai(\zeta) & \Ai(\omega^2 \zeta) \\
    \Ai'(\zeta) & \omega^2 \Ai'(\omega^2 \zeta) \end{array} \right)
    e^{- \frac{\pi i}{6} \sigma_3}
    \left(\begin{array}{cc} 1 & 0 \\ -1 & 1 \end{array} \right) &
    \mbox{for } \sigma < \arg \zeta < \pi, \\[10pt]
    \left(\begin{array}{cc}
    \Ai(\zeta) & -\omega^2 \Ai(\omega \zeta) \\
    \Ai'(\zeta) & -\Ai'(\omega \zeta) \end{array} \right)
    e^{-\frac{\pi i}{6} \sigma_3}
    \left(\begin{array}{cc} 1 & 0 \\ 1 & 1 \end{array} \right) &
    \mbox{for } -\pi < \arg \zeta < -\sigma, \\[10pt]
    \left(\begin{array}{cc}
    \Ai(\zeta) & -\omega^2 \Ai(\omega \zeta) \\
    \Ai'(\zeta) & -\Ai'(\omega \zeta) \end{array} \right)
    e^{- \frac{\pi i}{6} \sigma_3} &
    \mbox{for } -\sigma < \arg \zeta < 0,
    \end{array} \right.
\end{equation}
with $\omega = e^{2\pi i/3}$.
\end{proposition}
\begin{proof}
The proof is similar to \cite[p.1523--1525]{DKMVZ2}.
\end{proof}

\begin{remark}
An analysis of the proof in \cite{DKMVZ2} shows that the constant $C$ in
(\ref{eq59}) can be taken uniformly for $\sigma$ in a compact subset $K_{\sigma}$
of $(\pi/3, \pi)$, for $A$ in a compact subset $K_A$ of $(1, \infty)$,
and for $\delta$ uniformly in some interval $(\delta_0, \delta_1)$
depending on $K_{\sigma}$ and $K_A$.

A similar remark applies to the parametrix $\tilde{P}$ to be
constructed below.
\end{remark}

\subsection{Parametrix near $\bar\beta$}
A similar construction yields a parametrix $\tilde{P}$ in a neighborhood
$\tilde{U}_{\delta} = \{ z \mid |z-\bar\beta| < \delta\}$ that satisfies
the following Riemann--Hilbert problem.

\subsubsection*{Riemann--Hilbert problem for $\tilde{P}$:}
The problem is to determine, for a given $\delta > 0$ sufficiently small,
a $2\times 2$ matrix valued function
$\tilde{P}: \overline{\tilde{U}}_{\delta} \setminus \Sigma^T \to \mathbb C^{2\times 2}$
such that
\begin{enumerate}
\item[(a)] $\tilde{P}(z)$ is analytic for $z \in U_{\delta} \setminus \Sigma^T$,
and continuous on $\overline{\tilde{U}}_{\delta} \setminus \Sigma^T$,
\item[(b)] $\tilde{P}(z)$ possesses continuous boundary values for $z \in  \Sigma^T \cap \tilde{U}_{\delta}$,
denoted by $\tilde{P}_+(z)$ and $\tilde{P}_-(z)$, and
\begin{equation} \label{eq515}
    \tilde{P}_+(z) = \tilde{P}_-(z)
   \left(\begin{array}{cc} 0 & 1 \\ -1 & 0 \end{array} \right)
   \qquad \mbox{ for } z\in \Gamma \cap \tilde{U}_{\delta},
\end{equation}
\begin{equation} \label{eq516}
    \tilde{P}_+(z) = \tilde{P}_-(z) \left(\begin{array}{cc} 1 & 0 \\
    e^{2n \tilde{\phi}(z)} & 1 \end{array} \right)
    \qquad \mbox{ for } z \in (\Sigma_3 \cup \Sigma_4) \cap \tilde{U}_{\delta},
\end{equation}
and
\begin{equation} \label{eq517}
    \tilde{P}_+(z) = \tilde{P}_-(z) \left( \begin{array}{cc}
     1 & e^{-2n \tilde{\phi}(z)} \\
     0 & 1 \end{array} \right)
     \qquad \mbox{ for } z\in \Sigma_1 \cap \tilde{U}_{\delta}.
\end{equation}
Recall that the $\tilde{\phi}$-function is defined in (\ref{eq311}).
\item[(c)] There exists a constant $C > 0$ such that for every
$z \in \partial \tilde{U}_{\delta} \setminus \Sigma^T$,
\begin{equation} \label{eq518}
    \| \tilde{P}(z) N^{-1}(z)  - I \| \leq \frac{C}{n}.
\end{equation}
\end{enumerate}

There is a one-to-one analytic mapping $\zeta = \tilde{f}(z)$ such that
\begin{equation} \label{eq519}
    \frac{2}{3} (-\tilde{f}(z))^{3/2} = \tilde{\phi}(z).
\end{equation}
Then $\tilde{f}$ maps $\tilde{U}_{\delta}$ onto a convex neighborhood of $\zeta = 0$.

\begin{proposition}
The Riemann--Hilbert problem for $\tilde{P}$ is solved by
\begin{equation} \label{eq520}
    \tilde{P}(z) = \tilde{E}(z) \tilde{\Psi}^{\sigma}(n^{2/3} \tilde{f}(z))
    e^{n \tilde{\phi}(z) \sigma_3}
\end{equation}
with
\begin{equation} \label{eq521}
    \tilde{E}(z) = \sqrt{\pi} e^{\pi i/6}
    \left(\begin{array}{cc} 1 & 1 \\ i & -i \end{array} \right)
    \left(n^{1/6} (-\tilde{f}(z))^{1/4} a(z) \right)^{\sigma_3},
\end{equation}
and
\begin{equation} \label{eq522}
    \tilde{\Psi}^{\sigma}(\zeta) =
    \left(\begin{array}{cc} 1 & 0 \\ 0 & -1 \end{array} \right)
    \Psi^{\sigma}(-\zeta)
    \left(\begin{array}{cc} 1 & 0 \\ 0 & -1 \end{array} \right)
\end{equation}
\end{proposition}
\begin{proof}
This follows as in \cite[p.1527]{DKMVZ2}.
\end{proof}

\begin{remark}
As in \cite{DKMVZ2} there is a full asymptotic expansion of $PN^{-1}$
and $\tilde{P}N^{-1}$ in inverse powers of $n$. This expansion leads to
uniform asymptotic expansions for the generalized Laguerre polynomials.
In the present paper we will compute asymptotics for the polynomials
to first order, but we will not carry out the further computations
to determine the coefficients of a complete asymptotic expansion.
\end{remark}

\section{Final transformation $T \mapsto S$}
\setcounter{equation}{0}

Using $N$, $P$, and $\tilde{P}$, we define for every $n \in \mathbb N$,
\begin{equation} \label{eq61}
    S(z) = T(z) N(z)^{-1} \qquad
    \mbox{for } z \in \mathbb C \setminus
        (\Sigma^T \cup \overline{U}_{\delta} \cup \overline{\tilde{U}}_{\delta}),
\end{equation}
\begin{equation} \label{eq62}
    S(z) = T(z) P(z)^{-1} \qquad \mbox{for } z \in U_{\delta} \setminus \Sigma^T,
\end{equation}
\begin{equation} \label{eq63}
    S(z) = T(z) \tilde{P}(z)^{-1} \qquad \mbox{for } z \in \tilde{U}_{\delta} \setminus \Sigma^T.
\end{equation}

\begin{figure}[th!]
\centerline{\includegraphics*[width=10.8cm,height=8.4cm]{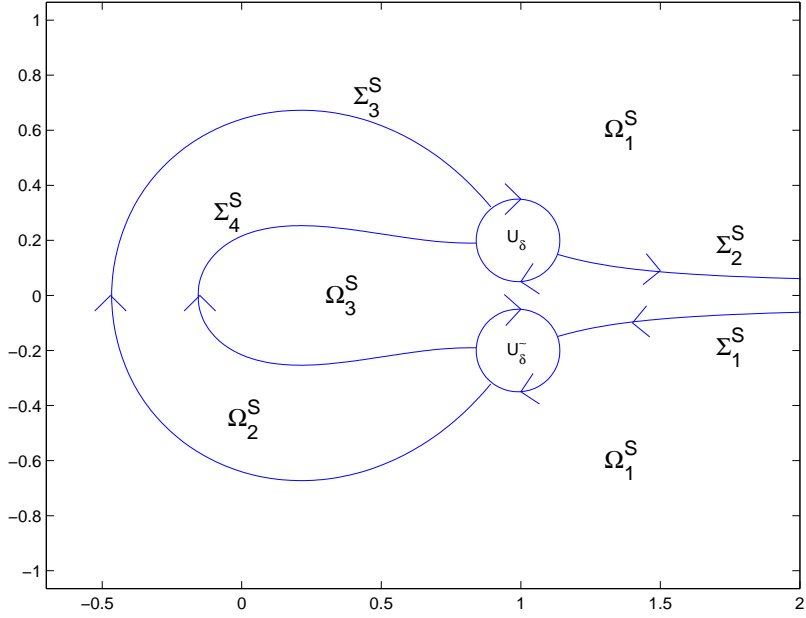}}
\caption{Contour $\Sigma^S = \partial U_{\delta} \cup \partial \tilde{U}_{\delta}
\cup \bigcup_j \Sigma_j^S$ and the domains $U_{\delta}$, $\tilde{U}_{\delta}$, and $\Omega_j^S$,
$j=1,2,3$, for the Riemann--Hilbert problem for $S$, for the value $A = 1.01$.}
\label{fig:plotS}
\end{figure}

Then $S$ is defined and analytic on $\mathbb C \setminus \left(\Sigma^T \cup \partial U_{\delta}
\cup \partial \tilde{U}_{\delta} \right)$. However it follows from the construction that
$S$ has no jumps on $\Gamma$ and on $\Sigma^T \cap (U_{\delta} \cup \tilde{U}_{\delta})$.
Therefore $S$ has an analytic continuation to $\mathbb C \setminus \Sigma^S$
(also denoted by $S$),
where $\Sigma^S$ is the contour indicated in Figure \ref{fig:plotS} for $A = 1.01$.

Formally, the contours $\Sigma_j^S$ are given by $\Sigma_j^S = \Sigma_j \setminus (\overline{U}_{\delta} \cup \overline{\tilde{U}}_{\delta})$
for $j =1,2,3,4$, and the domains $\Omega_j^S$ are given by
$\Omega_1^S = \Omega_1 \setminus
(\overline{U}_{\delta} \cup \overline{\tilde{U}}_{\delta})$, $\Omega_2^S =
( \Omega_2 \cup \Gamma \cup \Omega_3) \setminus
     (\overline{U}_{\delta} \cup \overline{\tilde{U}}_{\delta})$, and
$\Omega_3^S = \Omega_4  \setminus (\overline{U}_{\delta} \cup \overline{\tilde{U}}_{\delta})$.

Then $S$ satisfies the following Riemann--Hilbert problem.

\subsubsection*{Riemann--Hilbert problem for $S$:}
The problem is to determine $S : \mathbb C \setminus \Sigma^S \to \mathbb C^{2\times 2}$
such that the following hold.
\begin{enumerate}
\item[(a)] $S(z)$ is analytic for $z \in \mathbb C \setminus \Sigma^S$,
\item[(b)] $S(z)$ possesses continuous boundary values for $z \in \Sigma^S$,
denoted by $S_+(z)$ and $S_-(z)$, and
\begin{equation} \label{eq64}
    S_+(z) = S_-(z) P(z) N(z)^{-1}
    \qquad \mbox{for } z \in \partial U_{\delta},
\end{equation}
\begin{equation} \label{eq65}
    S_+(z) = S_-(z) \tilde{P}(z) N(z)^{-1}
    \qquad \mbox{for } z \in \partial \tilde{U}_{\delta},
\end{equation}
\begin{equation} \label{eq66}
    S_+(z) = S_-(z) N(z) \left(\begin{array}{cc}
        1 & 0 \\ e^{2n\phi(z)} & 1 \end{array} \right) N(z)^{-1}
    \qquad \mbox{for } z \in \Sigma_3^S \cup \Sigma_4^S,
\end{equation}
\begin{equation} \label{eq67}
    S_+(z) = S_-(z) N(z) \left(\begin{array}{cc}
        1 & e^{-2n\tilde{\phi}(z)} \\ 0 & 1 \end{array} \right) N(z)^{-1}
    \qquad \mbox{for } z \in \Sigma_1^S,
\end{equation}
\begin{equation} \label{eq68}
    S_+(z) = S_-(z) N(z) \left(\begin{array}{cc}
        1 & e^{-2n\phi(z)} \\ 0 & 1 \end{array} \right) N(z)^{-1}
    \qquad \mbox{for } z \in \Sigma_2^S,
\end{equation}
\item[(c)]  $S(z) = I + O\left(\frac{1}{z}\right)$ for $z \to \infty$.
\end{enumerate}

The jump matrices for $S$ are close to the identity matrix if $n$
is large. Indeed, by (\ref{eq59}) and (\ref{eq518}) we have
\[ \| P(z) N^{-1}(z) - I \| \leq \frac{C}{n} \qquad \mbox{for } z \in \partial U_{\delta}, \]
and
\[ \| \tilde{P}(z) N^{-1}(z) - I \| \leq \frac{C}{n}
    \qquad \mbox{for } z \in \partial \tilde{U}_{\delta}, \]
with a constant $C$ that is independent of $z$. The constant can also be chosen
independently of the value of $A$ for $A$ in a compact subset of $(1, \infty)$,
see Remark 5.4.

The jump matrices in (\ref{eq66})--(\ref{eq68}) are exponentially close to
the identity matrix.

\section{Asymptotic for polynomials $L_n^{(\alpha_n)}(nz)$}
\setcounter{equation}{0}

\subsection{Asymptotics for $S$}

In Sections 3--6 the values of $n$ and $\alpha$ were assumed to be fixed.
In order to study the asymptotics, we now let
$\alpha = \alpha_n$ depend on $n$ and we write $A_n = - \frac{\alpha_n}{n}$.
We assume that $A_n > 1$ for every $n$, and
\begin{equation} \label{eq71}
    \lim_{n \to \infty} - \frac{\alpha_n}{n} =
    \lim_{n \to \infty} A_n = A > 1.
\end{equation}
Because of the $n$-dependence of $A_n$, all of the notions and results introduced
in Sections 3--6 are $n$-dependent. For example, we have that
the curves $\Gamma$, $\Sigma$, $\Sigma^T$, and $\Sigma^S$ are all varying with $n$,
and so we denote them by $\Gamma_n$, $\Sigma_n$, $\Sigma^T_n$, and $\Sigma^S_n$.
Likewise, we have that the functions $R$, $g$, $\phi$, and $\tilde{\phi}$,
as well as all matrix-valued functions are $n$-dependent, and we also use a
subscript $n$ to denote their dependence on $n$.
If now we use $\Gamma$, $\Sigma$, $R$, $g$, etc., without
subscript $n$, then this refers to the limiting case.
Due to (\ref{eq71}) we have that the curves $\Gamma_n$ tend to the
limiting curve $\Gamma$, etc.

At the end of the previous section, we observed that the jump matrix
for $S_n$ is $I + O(1/n)$ uniformly on $\Sigma_n^S$ as $n \to \infty$.
In addition, the jump matrix converges to the identity matrix as $z \to \infty$
along the unbounded components of $\Sigma_n^S$ sufficiently fast, so that the
jump matrix is also close to $I$ in the $L^2$-sense. Since the contours $\Sigma_n^S$
are only slightly varying with $n$, we may follow arguments as those
given in \cite{BDJ} and \cite{Deift}, to conclude that
\begin{equation} \label{eq72}
    S_n(z) = I + O \left(\frac{1}{n}\right) \qquad \mbox{ as } n \to\infty
\end{equation}
uniformly for $z \in \mathbb C \setminus \Sigma_n^S$. Briefly, the Riemann--Hilbert
problem for $S_n$ is equivalent to a system of singular integral equations.
This system is of the form $\gamma(\zeta) - F(\gamma)(\zeta)= g(\zeta)$, where
$\gamma$ is a matrix valued function defined on the contour $\Sigma_n^S$, and
where the singular integral operator $F(\gamma)(\zeta)$ has operator norm
$O(\frac{1}{n})$. The system of singular integral equations can therefore
be solved by Neumann series in powers of $\frac{1}{n}$, and this leads to
a series expansion for the solution of the Riemann--Hilbert problem.

Using more precise information on the jump matrix for $S_n$ as indicated
in Remark 5.6, and with additional assumptions on the limit
(\ref{eq71}), one is able to obtain a full asymptotic expansion for
$S_n(z)$ in powers of $1/n$. This in return would give a full asymptotic
expansion for the generalized Laguerre polynomials.

\subsection{Strong asymptotics for generalized Laguerre polynomials}

Unraveling  the steps $Y_n \mapsto U_n \mapsto T_n \mapsto S_n$
and using (\ref{eq72}), we obtain strong asymptotics for $Y_n$ in all
regions of the complex plane. In particular we are interested in the
(1,1) entry of $Y_n$, since this is the monic generalized Laguerre
polynomial.

\begin{theorem}
Suppose for each $n \in \mathbb N$, we have a parameter $\alpha_n < -n$
such that {\rm (\ref{eq71})} holds where $A_n = - \frac{\alpha_n}{n}$.
Then we have the following asymptotic results for the generalized
Laguerre polynomials $L_n^{(\alpha_n)}(nz)$ as $n \to \infty$.
\begin{enumerate}
\item[\rm (a)] {\bf (asymptotics  away from $\Gamma$)} \\
Uniformly for $z$ in compact subsets of $\mathbb C \setminus \Gamma$,
we have as $n \to \infty$,
\begin{equation} \label{eq73}
    L_n^{(\alpha_n)}(nz) = \frac{(-n)^n}{n!} e^{ng_n(z)}
        \left( \frac{1 + R_n'(z)}{2} \right)^{1/2} \left(1 + O \left(\frac{1}{n}\right) \right).
\end{equation}
\item[\rm (b)] {\bf (asymptotics on $+$-side of $\Gamma_n$, away
from endpoints)} \\
Uniformly for $z$ on the $+$-side of $\Gamma_n$ away from $\beta$ and $\bar\beta$,
we have as $n \to \infty$,
\begin{eqnarray} \nonumber
    L_n^{(\alpha_n)}(nz) &  = & \frac{(-n)^n}{n!} e^{ng_n(z)}
       \left( \frac{1 + R_n'(z)}{2} \right)^{1/2} \\
       & & \times
       \left[ 1 - \left(\frac{1-R_n'(z)}{1+R_n'(z)} \right)^{1/2} e^{2n\phi_n(z)}
        + O \left(\frac{1}{n}\right) \right].
    \label{eq74}
\end{eqnarray}
\item[\rm (c)] {\bf (asymptotics on $-$-side of $\Gamma_n$, away from endpoints)} \\
Uniformly for $z$ on the $-$-side of $\Gamma_n$ away from $\beta$ and $\bar\beta$,
we have as $n \to \infty$,
\begin{eqnarray} \nonumber
    L_n^{(\alpha_n)}(nz) & = & \frac{(-n)^n}{n!} e^{ng_n(z)}
       \left( \frac{1 + R_n'(z)}{2} \right)^{1/2} \\
       & & \times
       \left[ 1 + \left(\frac{1-R_n'(z)}{1+R_n'(z)} \right)^{1/2} e^{2n\phi_n(z)}
        + O \left(\frac{1}{n}\right) \right].
    \label{eq75}
\end{eqnarray}
\item[\rm (d)] {\bf (asymptotics near $\beta$)} \\
Uniformly for $z$ in a (small) neighborhood of $\beta$, we have as $n \to \infty$,
\begin{eqnarray} \nonumber
\lefteqn{L_n^{(\alpha_n)}(nz)  = \frac{(-n)^n}{n!}
    \exp\left( \frac{n}{2} (A_n \log z + z  + \ell)\right) \sqrt{\pi} } \\
    & & \nonumber \times \,
    \left[ \left(\frac{z-\bar\beta_n}{z-\beta_n}\right)^{1/4} \left( n^{2/3} f_n(z) \right)^{1/4}
        \Ai(n^{2/3} f_n(z)) \left(1 + O\left(\frac{1}{n}\right) \right) \right. \\
    & & \qquad \label{eq76}
    \left. - \left(\frac{z-\beta_n}{z-\bar\beta_n} \right)^{1/4}
        \left( n^{2/3} f_n(z) \right)^{-1/4}
        \Ai'(n^{2/3} f_n(z)) \left(1 + O\left(\frac{1}{n}\right) \right) \right].
\end{eqnarray}
\end{enumerate}
\end{theorem}

\begin{proof}
(a) Let $K$ be a compact subset of $\mathbb C \setminus \Gamma$. Then, for
$n$ large enough, say $n \geq n_0$, we have $K \subset \mathbb C \setminus \Gamma_n$.
Let $n \geq n_0$. While choosing the contours $(\Sigma_3)_n$ and $(\Sigma_4)_n$
in Section 4, we then take care that $(\Omega_2)_n$
and $(\Omega_3)_n$ do not meet the compact $K$, see Figure \ref{fig:plotT}.
We also choose $\delta$ so that the disks of radius $\delta$ around $\beta_n$
and $\bar\beta_n$ are disjoint from $K$.
Then $K \subset (\Omega_1^S)_n \cup (\Omega_3^S)_n$, see Figure \ref{fig:plotS}.
For $z \in K$, we then  have by (\ref{eq41}), (\ref{eq48}), and (\ref{eq61})
\begin{eqnarray} \nonumber
    Y_n(z) & = & e^{n(\ell_n/2)\sigma_3} U_n(z) e^{ng_n(z)\sigma_3}
        e^{-n(\ell_n/2)\sigma_3} \\[10pt]
    \nonumber
    & = & e^{n(\ell_n/2)\sigma_3} T_n(z) e^{ng_n(z)\sigma_3} e^{-n(\ell_n/2)\sigma_3} \\[10pt]
    & = & e^{n(\ell_n/2)\sigma_3} S_n(z)N_n(z) e^{ng_n(z)\sigma_3} e^{-n(\ell_n/2)\sigma_3}.
    \label{eq77}
\end{eqnarray}
Using (\ref{eq28}), (\ref{eq55}), and (\ref{eq72}), and the fact that
$|(N_n)_{11}(z)|$ is bounded away from zero on $K$,  we obtain (\ref{eq73})
from (\ref{eq77}).

\medskip
(b) For $z \in  (\Omega_2)_n \setminus ((U_{\delta})_n \cup (\tilde{U}_{\delta})_n)$,
we have by (\ref{eq41}) and (\ref{eq49})
\begin{equation}
(Y_n)_{11}(z) = (U_n)_{11}(z) e^{ng_n(z)}
    =   \left[ (T_n)_{11}(z) + (T_n)_{12}(z) e^{2n\phi_n(z)} \right]
    e^{ng_n(z)} \label{eq78}
\end{equation}
Since $T_n = S_nN_n$ by (\ref{eq61}), $S_n = I + O(1/n)$, and the
entries of $N_n$ are uniformly bounded away from zero in the region
under consideration, we then get
\begin{equation}
(Y_n)_{11}(z) = \label{eq79}
        \left[ (N_n)_{11}(z) + (N_n)_{12}(z) e^{2n\phi_n(z)} + O(1/n) \right]
    e^{ng_n(z)},
\end{equation}
uniformly  for $z \in (\Omega_2)_n \setminus ((U_{\delta})_n \cup (\tilde{U}_{\delta})_n)$.
This proves (\ref{eq74}).

\medskip
(c) This is proved similarly as part (b). The only difference is that
for $z \in (\Omega_3)_n \setminus ((U_{\delta})_n \cup (\tilde{U}_{\delta})_n)$
we use (\ref{eq410}) instead of (\ref{eq49}).
This leads to
\[ (Y_n)_{11}(z) = \left[ (T_n)_{11}(z) - (T_n)_{12}(z) e^{2n\phi_n(z)} \right]
    e^{ng_n(z)} \]
instead of (\ref{eq78}). The rest of the proof is the same.

\medskip
(d) The proof is as in \cite[p.1539]{DKMVZ2}.
Let $U_{\epsilon}$ be an $\epsilon$-neighborhood of $\beta$ with $\epsilon < \delta$.
Then $U_{\epsilon} \subset (U_{\delta})_n$ for all large enough $n$, say $n \geq n_1$.
Let $z \in U_{\epsilon}$ and $n \geq n_1$.
For $z \in (\Omega_1)_n$, we find by (\ref{eq41}), (\ref{eq48}),
and (\ref{eq62}) that
\begin{eqnarray} \nonumber
    (Y_n)_{11}(z) & = &  (U_n)_{11}(z) e^{ng_n(z)} \, = \, (T_n)_{11}(z) e^{ng_n(z)} \\[10pt]
    \label{eq710}
    & = & \left[(S_n)_{11}(z) (P_n)_{11}(z) + (S_n)_{12}(z) (P_n)_{21}(z) \right] e^{ng_n(z)}.
\end{eqnarray}
By Proposition 5.3 we have
\[ P_n(z) = E_n(z) \Psi^{\sigma}(\zeta) e^{n\phi_n(z)\sigma_3} \]
where $\zeta = n^{2/3} f_n(z)$. For $z \in (\Omega_1)_n$, we have
that $\zeta$ belongs to the sector $0 < \arg \zeta < \sigma$,
so that the first formula for $\Psi^{\sigma}(\zeta)$
in (\ref{eq514}) applies. It follows that
\begin{equation} \label{eq711}
    \left(\begin{array}{c} (P_n)_{11}(z) \\ (P_n)_{12} (z) \end{array} \right)
    = E_n(z) \left(\begin{array}{c} \Ai(\zeta) \\ \Ai'(\zeta) \end{array} \right)
    e^{- \frac{\pi i}{6}} e^{n\phi_n(z)}.
\end{equation}
Using the definition (\ref{eq513}) of $E_n(z)$, we obtain from (\ref{eq711})
\[  \left(\begin{array}{c} (P_n)_{11}(z) \\ (P_n)_{12}(z) \end{array} \right)
     = \sqrt{\pi} \left(\begin{array}{cc} 1 & -1 \\ -i & -i \end{array} \right)
     \left(\frac{\zeta^{1/4}}{a_n(z)}\right)
    \left(\begin{array}{c} \Ai(\zeta) \\ \Ai'(\zeta) \end{array} \right)
   e^{n\phi_n(z)}.
 \]
Then
\begin{eqnarray}  \nonumber
\lefteqn{ (S_n)_{11}(z) (P_n)_{11}(z) + (S_n)_{12}(z) (P_n)_{21}(z) } \\
    && = \, \sqrt{\pi}  \left( \begin{array}{cc}
    1 + O(\frac{1}{n}) & -1 + O(\frac{1}{n}) \end{array} \right)
    \left(\begin{array}{c} \frac{\zeta^{1/4}}{a_n(z)} \Ai(\zeta) \\[10pt]
        \frac{a_n(z)}{\zeta^{1/4}} \Ai'(\zeta) \end{array} \right)
    e^{n\phi_n(z)}.\label{eq712}
\end{eqnarray}
Combining (\ref{eq710}) and (\ref{eq712}) with $\zeta = n^{2/3} f_n(z)$,
we obtain an expression
for $(Y_n)_{11}(z)$, which leads to (\ref{eq76}) in view of (\ref{eq24}) and (\ref{eq313}).

Next, for $z$ in the other regions $(\Omega_j)_n$, $j=2,3,4$,
similar calculations lead to the same expression for $(Y_{11})_n(z)$.
Hence (\ref{eq76}) holds uniformly for $z \in U_{\epsilon}$.

This completes the proof of Theorem 7.1.
\end{proof}

\begin{remark}
Since $\alpha_n$ is real, we have
$L_n^{(\alpha_n)}(n\bar{z}) = \overline{L_n^{(\alpha_n)}(nz)}$, and so
(\ref{eq76}) also describes the asymptotic behavior near $\bar\beta$.
\end{remark}

\subsection{Asymptotics for zeros}

From Theorem 7.1 we deduce the following results concerning the zeros
of the generalized Laguerre polynomials $L_n^{(\alpha_n)}(nz)$.

\begin{corollary}
Suppose we are in the same situation as in Theorem 7.1.
\begin{enumerate}
\item[\rm (a)] {\bf (All zeros tend to $\Gamma$)}
For every neighborhood $\Omega$ of $\Gamma$, there is $n_0$
such that for every $n \geq n_0$, all zeros of $L_n^{(\alpha_n)}(nz)$ are in $\Omega$.
\item[\rm (b)] {\bf (Zeros are on the $-$-side)}
For every $\delta > 0$, there is $n_0$ such that for every $n \geq n_0$,
there are no zeros of $L_n^{(\alpha_n)}(nz)$ in the region
$(\Omega_2)_n \setminus ((U_{\delta}) \cup (\tilde{U}_{\delta}))$.
\end{enumerate}
\end{corollary}

\begin{proof}
(a) This is immediate from the asymptotic formula (\ref{eq73}) since
$1+R'_n(z) \neq 0$.

\medskip
(b) Suppose $z$ is a zero of $L_n^{(\alpha_n)}(nz)$ lying in
$(\Omega_2)_n \setminus ((U_{\delta}) \cup (\tilde{U}_{\delta}))$,
that is, on the $+$-side of $\Gamma_n$. Then we have by (\ref{eq74})
\begin{equation} \label{eq713}
    \left(\frac{1-R'_n(z)}{1+R'_n(z)}\right)^{1/2} e^{2n \phi_n(z)} =
    1 + O\left(\frac{1}{n}\right).
\end{equation}
We know that $\Re \phi_n(z) \leq 0$. In order to obtain a contradiction
it is thus enough to enough to prove that
\begin{equation} \label{eq714}
    | 1 - R'_n(z)| < |1 + R'_n(z)|.
\end{equation}

Equality holds in (\ref{eq714}) if and only if
$R'_n(z)$ is purely imaginary, so that $(R'_n)^2(z)$ is real and negative.
From the definition (\ref{eq32}) of $R_n$, we get
\[ (R'_n)^2(z) = \frac{(z - \frac{\beta_n + \bar\beta_n}{2} )^2}{(z-\beta_n)(z-\bar\beta_n)}
    = \frac{(z - \Re \beta_n)^2}{(z-\Re \beta_n)^2 + (\Im \beta_n)^2}. \]
This is real and negative if and only if $z$ belongs to the vertical segment
connecting $\bar\beta_n$ and $\beta_n$. Consequently, this is the set where
$|1 - R_n'(z)| = |1 + R_n'(z)|$.
The vertical segment and $\Gamma$ form a closed contour, and it may be checked
that (\ref{eq714}) holds for $z$ outside of this contour, which
includes the region $(\Omega_2)_n \setminus ((U_{\delta}) \cup (\tilde{U}_{\delta}))$.

So we have a contradiction, and it follows that there are no zeros in
$(\Omega_2)_n \setminus ((U_{\delta}) \cup (\tilde{U}_{\delta}))$
for $n$ large enough.
\end{proof}

\begin{remark}
From the uniform asymptotics (\ref{eq76})  in a neighborhood
of $\beta$ it is possible to derive asymptotics for the extreme zeros
of $L_n^{(\alpha_n)}(nz)$. For example, it follows that for fixed $k \in \mathbb N$
and $n \geq k$, there is a zero $z_{k,n}$ of $L_n^{(\alpha_n)}(nz)$
such that
\begin{equation} \label{eq715}
    z_{k,n} = \beta_n - \frac{\iota_k}{f_n'(\beta_n) n^{2/3}} + O\left(\frac{1}{n} \right),
\end{equation}
where $-\iota_k$ is the $k$th largest (negative) zero of the Airy function
$\Ai(x)$.
\end{remark}

\noindent
\textit{A.B.J. Kuijlaars}
\hfill
{\footnotesize
 \texttt{arno@wis.kuleuven.ac.be}\\
 Department of Mathematics, Katholieke Universiteit Leuven,
 Celestijnenlaan 200 B, 3001 Leuven, Belgium
}

\noindent
\textit{K.T-R McLaughlin}
\hfill
{\footnotesize
 \texttt{mcl@amath.unc.edu}\\
 Department of Mathematics, University of North Carolina at Chapel Hill,
 Chapel Hill, NC 27599, USA\\
 and}
\hfill
{\footnotesize
 \texttt{mcl@math.arizona.edu}\\
 Department of Mathematics, University of Arizona, Tucson, AZ 85721, USA
}


\begin{thebibliography}{99}
\bibitem{AS}
    M. Abramowitz and I.A. Stegun,
    Handbook of mathematical functions,
    Dover Publications, New York, 1966.
\bibitem{BDJ}
    J. Baik, P. Deift, and K. Johansson,
    On the distribution of the length of the longest increasing
    subsequence of random permutations,
    J. Amer. Math. Soc. 12 (1999), 1119--1178.
\bibitem{BI}
    P. Bleher and A. Its,
    Semiclassical asymptotics of orthogonal polynomials,
    Riemann--Hilbert problem, and universality in the matrix model,
    Ann. Math. 150 (1999), 185--266.
\bibitem{dBSV}
    M. de Bruin, E.B. Saff, and R.S. Varga,
    On the zeros of generalized Bessel polynomials, I and II,
    Indag. Math. 43 (1981), 1--25.
\bibitem{Carp}
    A.J. Carpenter,
    Asymptotics for the zeros of the generalized Bessel polynomials,
    Numer. Math. 62 (1992), 465--482.
\bibitem{Deift}
    P. Deift,
    Orthogonal polynomials and random matrices: a Riemann--Hilbert approach,
    Courant Lecture Notes 3, Courant Institute 1999.
\bibitem{Deift2}
    P. Deift,
    Integrable systems and combinatorial theory,
    Notices Amer. Math. Soc. 47 (2000), 631--640.
\bibitem{DIZ}
    P. Deift, A. Its, and X. Zhou,
    A Riemann--Hilbert approach to asymptotic problems arising in the theory of
    random matrix models, and also in the theory of integrable statistical mechanics,
    Ann. Math. 146 (1997), 149--235.
\bibitem{DKMVZ1}
    P. Deift, T. Kriecherbauer, K.T-R McLaughlin, S. Venakides, and X. Zhou,
    Uniform asymptotics for polynomials orthogonal with respect to varying exponential
    weights and applications to universality questions in random matrix theory,
    Comm. Pure Appl. Math. 52 (1999), 1335--1425.
\bibitem{DKMVZ2}
    P. Deift, T. Kriecherbauer, K.T-R McLaughlin, S. Venakides, and X. Zhou,
    Strong asymptotics of orthogonal polynomials with respect to exponential weights,
    Comm. Pure Appl. Math. 52 (1999), 1491--1552.
\bibitem{DVZ}
    P. Deift, S. Venakides, and X. Zhou,
    New results in small dispersion KdV by an extension of the steepest
    descent method for Riemann--Hilbert problems.
    Internat. Math. Res. Notices 1997, no. 6, (1997), 286--299.
\bibitem{DZ1}
    P. Deift and X. Zhou,
    A steepest descent method for oscillatory Riemann--Hilbert problems,
    asymptotics for the mKDV equation,
    Ann. Math. 137 (1993), 295--368.
\bibitem{DZ2}
    P. Deift and X. Zhou,
    Asymptotics for the Painlev\'e II equation. Comm. Pure Appl. Math. 48 (1995), 277--337.
\bibitem{Dunster}
    T.M. Dunster,
    Uniform asymptotic expansions for the reverse generalized Bessel
    polynomials, and related functions,
    SIAM J. Math. Anal. 32 (2001), 987--1013.
\bibitem{FIK}
    A.S. Fokas, A.R. Its, and A.V. Kitaev,
    The isomonodromy approach to matrix models in 2D quantum gravity,
    Commun. Math. Phys. 147 (1992), 395--430.
\bibitem{GR}
    A.A. Gonchar and E.A. Rakhmanov,
    Equilibrium distributions and the rate of rational approximation of
    analytic functions,  Mat. Sbornik 134 (1987), 306--352.
    English translation in Math. USSR-Sbornik 62 (1989), 305--348.
\bibitem{Gros}
    E. Grosswald,
    Bessel polynomials, Lecture Notes Math. 698, Springer, New York, 1978.
\bibitem{KMM}
    S. Kamvissis, K.T-R McLaughlin, and P. Miller,
    Semiclassical soliton ensembles for the focusing nonlinear Schr\"odinger equation,
    manuscript.
\bibitem{KM}
    T. Kriecherbauer and K.T-R McLaughlin,
    Strong asymptotics of polynomials orthogonal with respect to Freud weights.
    Internat. Math. Res. Notices 1999, no. 6, (1999), 299--333.
\bibitem{MMO1}
    A. Mart\'{\i}nez-Finkelshtein, P. Mart\'{\i}nez-Gonz\'alez, and R. Orive,
    On asymptotic zero distribution of Laguerre and generalized Bessel polynomials
    with varying parameters,
    J. Comput. Appl. Math. 127 (2001), 255-266.
\bibitem{Olver1}
    F.W.J. Olver,
    The asymptotic expansion of Bessel functions of large order,
    Philos. Trans. Roy. Soc. London Ser. A 247 (1954), 328--368.
\bibitem{Olver2}
    F.W.J. Olver,
    Asymptotics and special functions,
    Academic Press, New York, 1974.
    Reprinted by AK Peters, Wellesley, 1997.
\bibitem{Perron}
    O. Perron, Die Lehre von den Kettenbr\"uchen, 3rd ed.
    Teubner, Stuttgart, 1957.
\bibitem{PV}
    I.E. Pritsker and R.S. Varga,
    The Szeg\H{o} curve, zero distribution and weighted approximation,
    Trans. Amer. Math. Soc. 349 (1997), 4085--4105.
\bibitem{SV1}
    E.B. Saff and R.S. Varga,
    On the zeros and poles of Pad\'e approximants to $e^z$,
    Numer. Math. 25 (1975), 1--14.
\bibitem{SV2}
    E.B. Saff and R.S. Varga,
    On the zeros and poles of Pad\'e approximants to $e^z$. III,
    Numer. Math. 30 (1978), 241--266.
\bibitem{Stahl}
    H. Stahl,
    Orthogonal polynomials with complex valued weight function, I and II,
    Constr. Approx. 2 (1986), 225--240, 241--251.
\bibitem{Strebel}
    K. Strebel,
    Quadratic differentials, Springer-Verlag, Berlin, 1984.
\bibitem{Szego1}
    G. Szeg\H{o},
    \"Uber eine Eigenschaft der Exponentialreihe,
    Sitzungsber. Berl. Math. Ges.  23 (1924), 50--64.
\bibitem{Szego}
    G. Szeg\H{o},
    Orthogonal polynomials, AMS Colloquium Publications 23,
    Amer. Math. Soc., Providence RI, 1939.
\bibitem{WZ}
    R. Wong and J.-M. Zhang,
    Asymptotic expansions of the generalized Bessel polynomials,
    J. Comput. Appl. Math. 85 (1997), 87--112.
\end{thebibliography}
\end{document}